# Existence and spatial limit theorems for lattice and continuum particle systems[*]

**Mathew D. Penrose**

*Department of Mathematical Sciences*
*University of Bath*
*Bath BA2 7AY*
*England*
*e-mail:* m.d.penrose@bath.ac.uk

**Abstract:** We give a general existence result for interacting particle systems with local interactions and bounded jump rates but noncompact state space at each site. We allow for jump events at a site that affect the state of its neighbours. We give a law of large numbers and functional central limit theorem for additive set functions taken over an increasing family of sub-cubes of $Z^d$. We discuss application to marked spatial point processes with births, deaths and jumps of particles, in particular examples such as continuum and lattice ballistic deposition and a sequential model for random loose sphere packing.

## Contents

---

[*]This is an original survey paper

## 1. Introduction

An *interacting particle system* may be informally described as 'a Markov process consisting of countably many pure-jump processes that interact by modifying each other's transition rates' ([24], p. 641). It is often assumed, as in [24] and most of [40], that these pure-jump processes live in a finite or at least compact state space; however, in applications this is not always the case. In the present work (in Section 2) we give a general construction for interacting particle systems where we assume only that the constituent pure-jump processes live in a separable complete metric space that is not required to be compact. We allow jumps occurring at a given site to affect the states (not just the jump rates) at neighbouring sites.

As well as a general existence result, we give general results on spatial limit theory. Typically, the interacting jump processes are indexed by the integer lattice $\mathbb{Z}^d$, with local interactions. This infinite system may be viewed as a limit of a sequence of finite systems, indexed by large finite 'windows' $W_n \subset \mathbb{Z}^d$. Averaging a (time-dependent) quantity of interest over the finite system $W_n$, one might expect to see it satisfy a law of large numbers and (functional) central limit theorem as $W_n \uparrow \mathbb{Z}^d$, and we give general results of this type in Section 3.

lattice and continuum models of sequential particle deposition. These have been much studied in the physical sciences literature [5; 21; 53; 61], and more recently in the mathematical literature [2; 6; 45; 46; 47; 49; 50; 51; 57]. In these, particles fall at random from above towards a surface (or *substrate*) and deposit themselves onto the surface (or fail to do so), in a manner which depends on the existing configuration of particles on the surface. According to the details of the model, the aggregation of deposited particles may be monolayer or multilayer, and the locations of deposited particles may be restricted to a lattice or may be continuous variables.

The results of the present paper extend known existence and limit results to functional central limit theorems, and to deposition models allowing for *displacement* whereby an incoming particle may influence the positions of existing particles. Models of 'random loose sphere packing' [9; 10; 14; 55] often feature displacement in some form. The framework given here encompasses both lattice and continuum models, both in monolayer and multilayer form. Our limit results apply to quantities such as the total number of deposited particles, or (in the case of multi-layer models) the total height of 'exposed' particles (i.e., those which lie at the surface of the accumulation of particles), and the average number of inter-particle contacts per unit area. Often in the physical sciences literature, such quantities are of interest, and their limiting values are assumed to exist without rigorous mathematical proof.

More generally, our results apply to any locally finite marked point process evolving jumpwise with locally determined jump rates and with locally bounded jump rates (here the word 'jump' could refer to births and deaths as well as movements of points). We refer to these processes as spatial *birth, death, migration and displacement* processes. They are described in general terms in Section 4.1 and special cases such as deposition models are discussed in the rest of Section 4. In this type of continuum model, the 'state' at a site in $\mathbb{Z}^d$ refers to the configuration of points in a patch of $\mathbb{R}^d$. Our results add to previous work in [26; 27; 33; 52] on the theory of evolving point process.

In recent complementary work, Garcia and Kurtz [25] define spatial birth-death processes in the continuum as solutions of so-called stochastic equations, and give criteria for existence and uniqueness; Qi [54] obtains functional central limit theorems for such processes (the work in [54] was carried out independently of that in this paper, which in turn was carried out independently of [25]). Unlike the present paper, [25] and [54] do not consder migration and displacement; on the other hand, their approach allows some relaxation of the assumption in the present paper of bounded range interactions.

We discuss some discrete-space examples and open problems in Section 5. Subsequent sections are devoted to proofs.

## 2. A general existence theorem

In this section we give our general result on existence and uniqueness of interacting particle systems. Before describing this, we briefly review some concepts from the theory of time-homogeneous Markov processes. Suppose $(E, \mathcal{E})$ is a measurable space. A *probability kernel* on $E$ is a mapping $\mu : E \times \mathcal{E} \to \mathbb{R}$ such that $\mu(\cdot, A)$ is measurable for each $A \in \mathcal{E}$ and $\mu(x, \cdot)$ is a probability measure on $E$ for each $x \in E$. A Markovian family of *transition distributions* on $E$ is a collection of probability kernels $(\mu_t, t \geq 0)$ with $\mu_0(x, \cdot) \equiv \delta_x$ for each $x \in E$ and with $\int_E \mu_t(y, A) \mu_s(x, dy) = \mu_{s+t}(x, A)$ for all $A \in \mathcal{E}$, $x, y \in E$, and $s, t \geq 0$. The induced (transition) *semigroup* of operators $(P_t, t \geq 0)$ is given by $P_t f(x) = \int_E f(y) \mu_t(x, dy)$, defined for all bounded measurable $f : E \to \mathbb{R}$. The *generator* of this semigroup is given by

$$Gf = \lim_{t \downarrow 0} t^{-1}(P_t f - f), \tag{2.1}$$

defined for all bounded measurable $f : E \mapsto \mathbb{R}$ for which the limit exists with uniform convergence. Note that in this paper we require uniform rather than pointwise convergence here for $Gf$ to be defined, as in e.g. [22] page 8 or [20] page 22.

Suppose that $E$ is a topological space and $\mathcal{E}$ is the Borel $\sigma$-field on $E$. Suppose $(\mu_t, t \geq 0)$ is a Markovian family of transition distributions on $E$, and $(\mathcal{G}_t, t \geq 0)$ is a filtration. A *Markov family of processes* in state-space $E$ with transition distributions $(\mu_t, t \geq 0)$ is a collection of stochastic processes $(\xi_t^x, t \geq 0)$, defined for each $x \in E$, which have right-continuous sample paths with left limits, are adapted to the filtration, and satisfy $P[\xi_0^x = x] = 1$ and $P[\xi_{s+t}^x \in A|\mathcal{G}_s] = \mu_t(\xi_s^x, A)$ almost surely, for each $(x, A) \in E \times \mathcal{E}$.

### 2.1. The general framework

Let $X$ be a Polish space (i.e., a separable, completely metrizable topological space). We consider Markov processes on the state space $X^V$, where $V$ is finite or countable (usually $\mathbb{Z}^d$ or finite subset thereof), endowed with the product topology. We assume that $V$ is the vertex set of a graph $\mathcal{V}$ with

$$D := \max\{\text{degree}(v) : v \in V\} < \infty. \tag{2.2}$$

If the state of the process at time $t$ is denoted $\xi_t$, then $\xi_t(v)$ is the state at site $v$ and lies in the space $X$; we shall refer to $X$ as the *local state space* and to $X^V$ as the *global state space*. We do *not* require $X$ to be finite or even compact; this feature puts us outside the scope of most of the existing interacting particle systems literature (though not all, as we discuss later). Many of the examples we consider are continuum systems where $X$ is required to be infinite.

Locally, the process evolves as a pure jump type Markov process. We assume *finite range*. That is, if a change occurs at site $v \in V$, then the instantaneous effect of this change is restricted to sites in a neighbourhood of $v$, and the jump rate at $v$ is determined by the state of the system in a neighbourhood of $v$. The neighbourhoods are determined by our assumed graph structure on $V$; for each $v \in V$, define $\mathcal{N}_v \subseteq V$ to consist of $v$ along with all adjacent vertices in the graph $\mathcal{V}$. The neighborhood $\mathcal{N}_v$ represents the set of all sites that can be instantaneously affected by changes at $v$ or which affect the jump rate at site $v$.

For each $v \in V$, assume we are given an operator $G_v^*$ which is the generator of a continuous-time pure jump type Markov process on the state space $X^{\mathcal{N}_v}$ with bounded rate function (for definitions see [38] page 238, or [24] page 635). Then for bounded measurable $f : X^{\mathcal{N}_v} \to \mathbb{R}$, and for $x \in X^{\mathcal{N}_v}$, we have (see [38], Proposition 19.2, or [24], eqn 31.11 or [22], page 376)

$$G_v^* f(x) = \int_{X^{\mathcal{N}_v}} (f(y) - f(x))\alpha_v^*(x, dy) \tag{2.3}$$

where the measure $\alpha_v^*(x, \cdot)$ is the jump rate kernel for the Markov process on $X^{\mathcal{N}_v}$ (the generator $G_v^*$ generates the process as it would proceed if changes took place only 'at $v$'). Let $\alpha_v^*(x)$ denote the total measure of $\alpha_v^*(x, \cdot)$, so that the next jump of the process with generator $G_v^*$, when in state $x$, occurs at rate $\alpha_v^*(x)$ and is governed by the measure $\alpha^*(x, \cdot)/\alpha_v^*(x)$ when it does occur. We assume throughout that rate functions are uniformly bounded, i.e. that

$$c_{\max} := \sup\{\alpha_v^*(x) : v \in V, x \in X^{\mathcal{N}_v}\} < \infty. \tag{2.4}$$

Let $\mathcal{C}$ be the class of bounded measurable functions $f : X^V \to \mathbb{R}$ which depend only on finitely many coordinates. Let $v \in V$. We define a generator $G_v$ of a jump process on $X^V$, as follows. Given $x \in X^V$, and $y \in X^{\mathcal{N}_v}$, let $x|v|y$ be the element of $X^V$ which agrees with $x$ outside $\mathcal{N}_v$, and agrees with $y$ inside $\mathcal{N}_v$. In other words, for $w \in V$ set

$$(x|v|y)(w) := \begin{cases} y(w) & \text{if } w \in \mathcal{N}_v \\ x(w) & \text{otherwise.} \end{cases}$$

For $f \in \mathcal{C}$, let the function $f_x^v : X^{\mathcal{N}_v} \to \mathbb{R}$ be given by

$$f_x^v(y) := f(x|v|y). \tag{2.5}$$

Let $x|_{\mathcal{N}_v}$ be the restriction of $x$ to $\mathcal{N}_v$ (an element of $X^{\mathcal{N}_v}$), and let $x|_{V \setminus \mathcal{N}_v}$ be the restriction of $x$ to $V \setminus \mathcal{N}_v$. Set

$$G_v(f)(x) = G_v^*(f_x^v(x|_{\mathcal{N}_v})). \tag{2.6}$$

Then $G_v$ is the generator of a jump process with jump kernel denoted $\alpha_v$, given by

$$\alpha_v(x, x|v|dy) = \alpha_v^*(x|_{\mathcal{N}_v}, dy);$$
$$\alpha_v(x, dz) = 0 \quad \text{if} \quad z|_{V \setminus \mathcal{N}_v} \neq x|_{V \setminus \mathcal{N}_v}.$$

Here, for $z \in X^V$ the notation $z|_{V \setminus \mathcal{N}_v} \neq x|_{V \setminus \mathcal{N}_v}$ means there exists $w \in V \setminus \mathcal{N}_v$ with $z(w) \neq x(w)$.

Given $W \subseteq V$, define the operator $G^W$ on $\mathcal{C}$ by

$$G^W = \sum_{z \in W} G_z. \tag{2.7}$$

We are interested only in cases where $W$ is finite or $W = V$. In the latter case, write simply $G$ for $G^V$.

If $W \subset V$ is finite, it is clear that $G^W$ is the generator of a jump process with jump kernel given by the sum over $v \in W$ of the jump kernels $\alpha_v$. Write $P_t^W$ for the associated semigroup.

### 2.2. Existence and uniqueness

Let the space $X^V$ be endowed with the product topology and Borel $\sigma$-field. Our first general result gives existence and uniqueness for a Markov transition semigroup on $X^V$ with generator $G$, and existence of a Markov family of processes in state space $X^V$ corresponding to this transition semigroup. This process has values in $D([0, \infty), X^V)$, and is called an *interacting particle system*, at least when $X$ is finite (see [24], Chapter 32). Given a sequence of sets $(W_m, m \geq 1)$ we write $\liminf_{m \to \infty}(W_m)$ for $\cup_{m=1}^\infty \cap_{n=m}^\infty W_n$.

**Theorem 2.1.** (Existence and uniqueness result) *Let $G$ be the generator defined by (2.7) with $W = V$. Assume (2.2) and (2.4) hold. Then there exists a unique Markovian family of transition distributions $(\mu_t, t \geq 0)$ on $X^V$, such that the associated semigroup $(P_t, t \geq 0)$ has generator $G$ on $\mathcal{C}$.*

*Moreover, for any sequence $(W_m)_{m \geq 1}$ of finite subsets of $V$, and any $t \geq 0$*

$$P_t f(x) = \lim_{m \to \infty} P_t^{W_m} f(x), \ \forall \ f \in \mathcal{C}, x \in X^V, \ \text{if } \liminf_{m \to \infty}(W_m) = V. \qquad (2.8)$$

*There also exists a filtration $(\mathcal{G}_t, t \geq 0)$ and a Markov family of processes $(\xi_\cdot^x, x \in X^V)$, adapted to $(\mathcal{G}_t, t \geq 0)$, with transition semigroup $(P_t, t \geq 0)$.*

Most previous existence results for interacting particle systems have been given for cases when the space $X$ is taken to be finite (for example [24; 39]), or at least compact ([32], and Chapter 1 of [40]). Unlike in these results, we do not assume $X$ is compact. If $X$ is not compact then $X^V$ is not even locally compact, which also renders inapplicable much of the general literature on Markov process theory. Compactness of $X$ would imply existence of at least one equilibrium distribution (see, e.g., [24]), but we do not here concern ourselves with equilibrium distributions.

There is some previous literature in the case of non-compact $X$. In Chapter 9 of Liggett [40], the space $X$ is $[0, \infty)$ and the interactions have a linear structure. Our setup allows general noncompact $X$ and does not require this linear structure, although conversely, Chapter 9 of [40] does not require the interactions have the local structure considered here. Basis [7; 8] gives a general existence result in non-compact spaces, and Chen ([12], Chapter 13) gives a similar result, and also a uniqueness result. However, the general approach, and the detailed set of conditions, in those works seem very different from here. In particular, both [7] and [12] require a 'ground state' in $X^V$ to be specified, and take each approximating jump process, with semigroup $(P_t^{W_n}, t \geq 0)$, to have its state outside $W_n$ to be given by the ground state (which can be viewed as a 'boundary condition'). With our approach, there is no need to specify a ground state (and the approximating jump processes are not needed). The examples considered here are mainly in continuous spaces while those considered in [7; 12] are mainly in discrete spaces.

The literature includes a number of papers on existence and other properties of particular interacting particle systems with non-compact $X$ and non-linear structure (in some cases with unbounded jump rates). For example, [41] gives a construction for several such processes and the method can be adapted to others with non-linear rates (see [1]). In these examples the state space is not the whole of $X^V$ but a space of functions in $X^V$ satisfying a summability condition. Other examples include [58] and [4].

Our proof for Theorem 2.1 loosely follows the graphical representation of Harris [31] for finite-range particle systems, even though that paper, and subsequent presentations such as those in [19; 28; 40], are also restricted to the case where $X$ is finite. In extending this argument to a more general class of sets $X$, we use a measure-theoretic result from Kallenberg [38] which is used there to express Markov chains as random dynamical systems. Using this result, any pure jump continuous-time Markov process with bounded jump rates in a Polish state space can be generated by a Poisson process in $(0, \infty) \times [0, 1]$. In proving Theorem 2.1, we shall show that the processes $\xi_t^x$ generated by $G$ can be realized in terms of a family of independent Poisson processes, indexed by $v \in V$.

## 3. Spatial limit theorems

### 3.1. Further assumptions and definitions

We use the following notation throughout. Given $d \in \mathbb{N}$ we let $\mathbf{0}$ denote the origin in $\mathbb{R}^d$. For $A \subset \mathbb{R}^d$ and $u \in \mathbb{R}^d$ we write $u + A$ for the set $\{u + u' : u' \in A\}$. For $r > 0$, let $B_r$ be the closed Euclidean ball of radius $r$ centred at the origin in $\mathbb{R}^d$.

For the results in this section, we make some further assumptions, in addition to those of Section 2. We assume that $V = \mathbb{Z}^d$, that some specified 'window' $\mathcal{N}$ around the origin (i.e. a finite symmetric set $\mathcal{N} \subset \mathbb{Z}^d$ with $\mathbf{0} \in \mathcal{N}$) is given, and that $\mathcal{N}_v$ is the set $v + \mathcal{N}$ for each $v \in \mathbb{Z}^d$. We also assume translation invariance of the jump rates associated with each lattice point. In other words, we assume we are given a jump rate kernel $\alpha^*(x, dy)$ on $X^{\mathcal{N}}$, and that for $v \in \mathbb{Z}^d$, the jump kernel $\alpha_v^*$ on $X^{\mathcal{N}_v}$ is the shifted version of $\alpha^*$ given by

$$\alpha_v^*(x; \Gamma) = \alpha^*(L_v(x), L_v(\Gamma)), \quad x \in X^{\mathcal{N}_v}, \ \Gamma \text{ Borel in } X^{\mathcal{N}_v}, \quad (3.1)$$

where $L_v : X^{\mathcal{N}_v} \to X^{\mathcal{N}}$ is given by

$$(L_v x)(w) = x(w + v), \quad x \in X^{\mathcal{N}_v}, w \in \mathcal{N}.$$

We assume the total measure of $\alpha^*(x, dy)$ is bounded, uniformly over $x \in X^{\mathcal{N}}$; this is enough to ensure that the condition (2.4) holds here.

Let $\mathcal{B}$ denote the collection of all non-empty finite subsets of $\mathbb{Z}^d$. Let $\nu$ be a probability measure on $X$. For $A \in \mathcal{B}$, let $(\xi_t^{A,\nu}, t \geq 0)$ be a realization of the $X^{\mathbb{Z}^d}$-valued pure jump Markov process with generator $G^A$ given by (2.7) and with initial value $\xi_0^{A,\nu} = (\xi_0^{A,\nu}(v), v \in \mathbb{Z}^d)$ whose entries $\xi_0^{A,\nu}(v)$ ($v \in \mathbb{Z}^d$) are independent and identically distributed states with common distribution $\nu$. Likewise, let $(\xi_t^{\nu}, t \geq 0)$ be a realization of the $X^{\mathbb{Z}^d}$-valued Markov process with generator $G$ and satisfying (2.8) (as given in Theorem 2.1) with initial value $\xi_0^{\nu}$ whose entries $\xi_0^{\nu}(v)$ ($v \in \mathbb{Z}^d$) are independent and identically distributed states with common distribution $\nu$.

Our limit theorems are given in terms of *stationary additive set functions* on $X^{\mathbb{Z}^d}$, which we define as follows. Given a measurable real-valued function $H : X^{\mathcal{N}} \to \mathbb{R}$, the corresponding stationary additive set function $S_H^A$, defined for each $A \in \mathcal{B}$, is a real-valued functional on the global state space $X^{\mathbb{Z}^d}$ obtained by summing $H$ over the states at (neighbourhoods of) sites in $A$, i.e., defined for $x \in X^V$ by

$$S_H^A(\xi) = \sum_{v \in A} H(L_v(x|_{\mathcal{N}_v})). \quad (3.2)$$

We shall give limit theorems, in the form of a law of large numbers (LLN) and a central limit theorem (CLT) as $A$ becomes large, for the random variable $S_H^A(\xi_t^{A,\nu})$. The idea here is that the evolution of the interacting particle system is restricted to the finite set $A$, and the function $H$ is also summed over (neighbourhoods of) states at sites in $A$. See Section 4 for some examples.

For $A \subset \mathbb{Z}^d$ we write $A^c$ for $\mathbb{Z}^d \setminus A$, and define the 'neighbourhood', 'exterior boundary', and 'interior' of $A$, respectively, by

$$\mathcal{N}_A := \cup_{v \in A} \mathcal{N}_v; \quad \partial_{\text{ext}} A := \mathcal{N}_A \setminus A; \quad A^o := \mathbb{Z}^d \setminus \mathcal{N}_{A^c}.$$

Also let $|A|$ denote the number of elements of $A$.

### 3.2. Limit theorems

For $\gamma > 0$ and $\tau \geq 0$, we shall consider $H$ and $\nu$ satisfying the moments condition

$$\sup_{A \in \mathcal{B}} E[|H(\xi_\tau^{A,\nu}|_{\mathcal{N}})|^\gamma] < \infty. \tag{3.3}$$

We consider $\mathcal{B}$-valued sequences $(A_n, n \geq 1)$ satisfying the conditions

$$|\partial_{\text{ext}} A_n|/|A_n| \to 0 \quad \text{as} \quad n \to \infty; \tag{3.4}$$

$$(\mathcal{N}_{A_n})^o = A_n \quad \forall n. \tag{3.5}$$

The conditions (3.4) and (3.5) hold, for example, if $A_n = ([-n, n] \cap \mathbb{Z})^d$.

We have the following LLN for functionals on interacting particle systems with finite range interactions.

**Theorem 3.1.** (General LLN) *Let $d$, $\mathcal{N}$ and $H$ be as in Section 3.1. Let $\tau > 0$. Let $\nu$ be a probability measure on $X$, and suppose that the moments condition (3.3) holds for some $\gamma > 1$. Then for any $\mathcal{B}$-valued sequence $(A_n, n \geq 0)$ satisfying (3.4) and (3.5), we have as $n \to \infty$ that*

$$|A_n|^{-1} S_H^{A_n}(\xi_\tau^{A_n,\nu}) \xrightarrow{L^1} EH(\xi_\tau^\nu|_{\mathcal{N}}). \tag{3.6}$$

We have the following CLT for functionals on interacting particle systems with finite range interactions, evaluated at a finite set of times.

**Theorem 3.2.** (General CLT) *Let $d$, $\mathcal{N}$ and $H$ be as in Section 3.1. Suppose $I \subseteq [0, \infty)$. Suppose that for some $\gamma > 2$, (3.3) holds for all $\tau \in I$. Then for all $s \in I, t \in I$, the sum*

$$\sigma(s,t) := \sum_{z \in \mathbb{Z}^d} \text{Cov}\left[H(\xi_s^\nu|_{\mathcal{N}}), H(L_z(\xi_t^\nu|_{\mathcal{N}_z}))\right] \tag{3.7}$$

*converges absolutely, and for any $\mathcal{B}$-valued sequence $(A_n)_{n \geq 1}$ satisfying (3.4) and (3.5), as $n \to \infty$ we have*

$$|A_n|^{-1} \text{Cov}(S_H^{A_n}(\xi_s^{A_n,\nu}), S_H^{A_n}(\xi_t^{A_n,\nu})) \to \sigma(s,t) \tag{3.8}$$

*and the finite-dimensional distributions of the process*

$$|A_n|^{-1/2}(S_H^{A_n}(\xi_t^{A_n,\nu}) - ES_H^{A_n}(\xi_t^{A_n,\nu})), \quad t \in I \tag{3.9}$$

*converge to those of a zero-mean Gaussian process with covariances given by $(\sigma(s,t), s,t \in I)$.*

When $I = [0, \tau]$ for some $\tau > 0$, it is possible to extend the convergence in distribution of finite dimensional distributions to a *functional* CLT giving weak convergence in $D[0, \tau]$ of the whole process $(S_H^{A_n}(\xi_t^{A_n, \nu}), 0 \leq t \leq \tau)$. Here $D[0, \tau]$ is the space of right-continuous functions on $[0, \tau]$ with left limits, equipped with the Skorohod topology (see [11]). One reason for an interest in such convergence is the study of weak convergence of hitting times (see [18]), although we do not pursue this further here. The functional CLT requires a stronger version of the moments condition (3.3), namely

$$\sup_{A \in \mathcal{B}} \sup_{0 \leq t \leq \tau} E[|H(\xi_t^{A,\nu}|_{\mathcal{N}})|^\gamma] < \infty. \tag{3.10}$$

**Theorem 3.3.** *Let $d$, $\mathcal{N}$ and $H$ be as in Section 3.1. Let $\tau > 0$. Suppose (3.10) holds for some $\gamma > 8$. Then with $I = [0, \tau]$, the process defined by (3.9) converges weakly in $D[0, \tau]$ to a zero-mean Gaussian process with covariance structure given by $\sigma(s, t), s, t \in [0, \tau]$. Also, the limiting process admits a version with continuous sample paths.*

### 3.3. Remarks

1. The limiting covariance in (3.8) does not depend on the choice of $(A_n)_{n \geq 1}$.
2. We do not in general rule out the possibility that the limiting covariance in (3.8) might be zero, although in most applications this will not be the case.
3. Theorem 3.2 can be extended to a CLT of measures, using results from [48]. Let $\mathcal{R}$ denote the class of Riemann measurable sets in $\mathbb{R}^d$ (see [48] for definitions), fix $A_0 \in \mathcal{R}$ with non-empty interior, and let $A_1, \ldots, A_k$ be Riemann measurable subsets of $A_0$. Let $(s_n, n \geq 1)$ be an increasing unbounded sequence of positive numbers, and let $A_{i,n} := (s_n A_i) \cap \mathbb{Z}^d$. Let $t_1, \ldots, t_k$ be positive numbers. Then if (3.3) holds for some $\gamma > 2$ and $\tau \geq \max(t_1, \ldots, t_n)$, we expect that

$$(s_n^{-d/2}(S_H^{A_{i,n}}(\xi_{t_i}^{A_{0,n}, \nu}) - E[S_H^{A_{i,n}}(\xi_{t_i}^{A_{0,n}, \nu})]), \quad 1 \leq i \leq k)$$

   converges, as $n \to \infty$, to a multivariate centred normal with covariances given by $\sigma(t_i, t_j) \times \text{Leb}(|A_i \cap A_j|)$ where $\text{Leb}(\cdot)$ denotes Lebesgue measure. It should be possible to prove this using theorem 2.1 of [48], but we have not written out the details.
4. An alternative approach to the CLT would be via Stein's method, as used by Penrose and Sudbury ([49], theorem 7) for certain specific particle systems with finite local state space. It may be possible to adapt that approach to give a normal approximation in the more general setting of Section 2 (i.e., without assuming $V$ is $\mathbb{Z}^d$ and with no stationarity assumptions), so long as we simplify $H$ to a function $\tilde{H}$ depending only on the local state (and not on neighbouring sites). That is, taking $\tilde{H}$ to be a function from $X$ to $\mathbb{R}$, in the general setting of Section 2, one might expect normal approximation for

$$S^A_{\tilde{H}}(\xi^{A,\nu}_t) := \sum_{v \in A} \tilde{H}(\xi^{A,\nu}_t(v)).$$

Again, we have not written out the details.

5. These results are related to those given in Penrose [47] and in Doukhan et al. [18] for interacting particle systems with finite local state spaces. For other approaches to central limit theorems for interacting particle systems, see Holley and Stroock [34]. See also De Masi and Presutti [16], Spohn [60], and Liggett ([40], page 39). These works are mainly concerned with cases where $X$ is finite, and where there exists a unique invariant measure, which we do *not* assume here.

6. The initial condition on the process $\xi^\nu_t$, namely independent and identically distributed states at time zero, is perhaps more restrictive than one would like. We believe that it is possible to give similar results to Theorems 3.2 and 3.3 in the more general case of a spatially stationary initial distribution satisfying a strong mixing condition, if instead of the process $\xi^A_t$ we consider $\xi_t$, thereby making everything spatially stationary. Some other ways in which one might hope to extend the current work are described in Section 5.2.

## 4. Spatial birth, death, migration and displacement processes

### 4.1. General description

Our setup is intended to be applicable to a large class of continuum particle systems. In these, the global state space $\mathcal{S}$ is the space of locally finite subsets of $\mathbb{R}^d$, $d \in \mathbb{N}$.

Let $R > 0, \lambda > 0$ be fixed parameters; $R$ denotes the range of interaction. Let $\mathcal{S}_R$ be the space of finite configurations of points in the ball $B_R$ centred at **0**. Given an element $\mathcal{X}$ of the state space $\mathcal{S}$, for $\mathbf{x} \in \mathbb{R}^d$ let $\mathcal{X}|_\mathbf{x}$ denote the set $(-\mathbf{x} + \mathcal{X}) \cap B_R$ (an element of $\mathcal{S}_R$). Consider the following possible transitions of the process when in state $\mathcal{X}$.

(i) *Immigration* ('birth'). Particles arrive as a space-time homogeneous Poisson process of rate $\lambda$. Incoming particles may be accepted or rejected. The probability of acceptance of an incoming particle at $\mathbf{x}$ is assumed to be a function of the point process $\mathcal{X}|_\mathbf{x}$. If accepted, the incoming particle is placed at $\mathbf{x}$ and added to the configuration.

(ii) *Death*. A particle at $\mathbf{x}$ disappears at rate $\delta(\mathcal{X}|_\mathbf{x})$ where for each $\mathcal{Y} \in \mathcal{S}_R$, $\delta(\mathcal{Y})$ is a nonnegative number.

(iii) *Migration*. A particle at $\mathbf{x}$ jumps at rate $\rho(\mathcal{X}|_\mathbf{x})$, and when it jumps its new location $\mathbf{x} + \mathbf{y}$ is distributed with distribution $\mu(\mathcal{X}|_\mathbf{x}; d\mathbf{y})$, where for each configuration $\mathcal{Y} \in \mathcal{S}_R$, $\rho(\mathcal{Y})$ is a nonnegative number and $\mu(\mathcal{Y}; \cdot)$ is a probability measure on $B_R$.

To these ingredients, we can add the following embellishments:

(a) *Marked point process.* Assume each particle carries a 'mark' in a measurable space $\mathcal{T}$ (the mark space). Incoming particles have random marks, assumed independent and identically distributed according to some probability distribution on $\mathcal{T}$. The probability of acceptance of an particle at $\mathbf{x}$ may depend on its mark and those of the existing points in $\mathcal{X}|_{\mathbf{x}}$. A 'jump' (migration) of a particle may entail a change of the value of its mark as well as (or instead of) a change of location in $\mathbb{R}^d$. In the marked case, we modify $\mathcal{S}$ to be the set of locally finite subsets of $\mathbb{R}^d \times \mathcal{T}$; however, we shall often refer to subsets of $\mathbb{R}^d \times \mathcal{T}$ more informally as 'marked subsets of $\mathbb{R}^d$'.

(b) *Synchronous updating.* Whenever any one of the above types of event (immigration, death, or migration) occurs, at a point $\mathbf{x} \in \mathbb{R}^d$, then this can cause some rearrangement, or change in the values of the marks, of the points within distance $R$ of $\mathbf{x}$ (but without any change to the configuration at a distance more than $R$ from $\mathbf{x}$). It is assumed that this rearrangement occurs simultaneously with the immigration/migration/death event, and is determined by the nature of that event and the configuration of (marked) points within distance $R$ of $\mathbf{x}$. In the case of an immigration event at $\mathbf{x}$, we allow the location and mark of the incoming particle itself to be affected by the existing points within distance $R$ of $\mathbf{x}$, so the particle may be placed at some point other than $\mathbf{x}$ but within distance $R$ of $\mathbf{x}$.

Continuum processes of this type fit into the particle system picture of Sections 2 and 3 as follows. Take the local state space $X$ to consist of all finite subsets of $[0,1)^d \times \mathcal{T}$ (or more formally, finite point measures on $[0,1)^d \times \mathcal{T}$). In the unmarked case we may think of $\mathcal{T}$ as having a single element. We partition $\mathbb{R}^d$ into cubes $C_v := v + [0,1)^d, v \in \mathbb{Z}^d$, and the global state space $\mathcal{S}$ is identified with $X^{\mathbb{Z}^d}$ by identifying a locally finite marked point set $\mathcal{X} \subset \mathbb{R}^d \times \mathcal{T}$ with an element $\eta_\mathcal{X}$ of $X^{\mathbb{Z}^d}$ obtained by setting

$$\eta_\mathcal{X}(v) = -v + (\mathcal{X} \cap (C_v \times \mathcal{T})), \quad v \in \mathbb{Z}^d,$$

where for $v \in \mathbb{Z}^d$ and $\mathcal{Y} \subset \mathbb{R}^d \times \mathcal{T}$ we set $v + \mathcal{Y} := \{(v + \mathbf{x}, \theta) : (\mathbf{x}, \theta) \in \mathcal{Y}\}$. More informally, the local state of the system at site $v \in \mathbb{Z}^d$ is the collection of marked points in $C_v$.

For $v \in \mathbb{Z}^d$, we take $\mathcal{N}_v$ to be the set of $w \in \mathbb{Z}^d$ such that the distance between the cubes $C_v$ and $C_v$ is at most $R$ (and take $\mathcal{N} := \mathcal{N}_0$). When an immigration, death or migration event occurs at a point $\mathbf{x} \in C_v$, with or without synchronous updating, this may affect the states of other vertices, but only those in $\mathcal{N}_v$. Therefore, events at vertex $v \in \mathbb{Z}^d$ affect only the states at vertices in $\mathcal{N}_v$, and this type of process fits into the setup of Section 2; in particular, (2.2) holds.

Moreover, rates are defined in a translation-invariant manner so that (3.1) holds.

With the above identification of $\mathcal{S}$ with $X^{\mathbb{Z}^d}$, a Markov process in $X^{\mathbb{Z}^d}$ corresponds to a Markov process in $\mathcal{S}$. Given the local interactions (i.e., the parameters of the immigration, death and migration rates), we can use Theorem 2.1 to deduce existence of a Markov family processes $(\mathcal{X}_t, t \geq 0)$ in the state-space $\mathcal{S}$ with these local interactions. For bounded Borel $\mathbf{A} \subseteq \mathbb{R}^d$, let $(\mathcal{X}_t^{\mathbf{A}}, t \geq 0)$ denote the corresponding process when births, deaths and migrations taking place outside $\mathbf{A}$ are suppressed.

A large class of stationary additive set functions on $\mathcal{S}$ is given as follows. Let $\phi : \mathcal{T} \times \mathcal{S}_R \to \mathbb{R}$ be measurable and for bounded Borel $\mathbf{A} \subset \mathbb{R}^d$, and $\mathcal{X} \in \mathcal{S}$ set

$$S_\phi^{\mathbf{A}}(\mathcal{X}) := \sum_{(\mathbf{x},\theta) \in \mathcal{X} \cap (\mathbf{A} \times \mathcal{T})} \phi(\theta, (-\mathbf{x} + \mathcal{X}) \cap (B_R \times \mathcal{T})).$$

The corresponding additive set function on $X^{\mathbb{Z}^d}$ is obtained by taking

$$H_\phi(\eta_\mathcal{X}|_\mathcal{N}) := \sum_{(\mathbf{x},\theta) \in \mathcal{X} \cap (C_{\mathbf{0}} \times \mathcal{T})} \phi(\theta, (-\mathbf{x} + \mathcal{X}) \cap (B_R \times \mathcal{T})),$$

which depends on $\eta_\mathcal{X}$ only through $\eta_\mathcal{X}|_\mathcal{N}$. Given $A \in \mathcal{B}$, let $\tilde{A}$ denote the set $\cup_{v \in A} C_v$. We have the identity

$$S_\phi^{\tilde{A}}(\mathcal{X}) = S_{H_\phi}^A(\eta_\mathcal{X}),$$

so if $H_\phi$ satisfies the moment condition (3.3) for appropriate $\gamma$, we can use Theorems 3.1, 3.2 and 3.3 to obtain corresponding limit theorems for $S_\phi^{\tilde{A}_n}(\mathcal{X}_t^{\tilde{A}_n})$.

The statements of results in Section 3 refer to processes $(\xi_t^{A,\nu}, t \geq 0)$, taking values in $X^{\mathbb{Z}^d}$, for finite $A \subset \mathbb{Z}^d$. In the $\mathcal{S}$-valued process $(\mathcal{X}_t^{\tilde{A},\nu})$ corresponding to $\xi_t^{A,\nu}$ the immigration, death and migration events occur only in the bounded region $\tilde{A} := \cup_{v \in A} C_v$, but within this region the description of rates of these events is as before; we reserve the notation $\xi_t$ for $X^{\mathbb{Z}^d}$-valued processes and $\mathcal{X}_t$ for processes in $\mathcal{S}$. Typically the physical system being modelled is in a bounded region, and Theorems 3.1, 3.2 and 3.3 provide information about the limit behaviour of the bounded-region process as this region becomes large. To apply these results, we need to specify $\nu$; that is, we need to assume the initial configurations of (marked) points in each cube $C_z$, $z \in A$ are independent and identically distributed with some specified distribution $\nu$. For example, one could assume an initially empty configuration, i.e. assume $\nu$ is concentrated on the configuration in $C_{\mathbf{0}}$ with no points.

To apply our general results, we need (2.4) to hold, i.e. we need a bound on the total jump rate at each site $v \in \mathbb{Z}^d$. With our description, the rate of immigration events per unit volume is bounded by $\lambda$; however, if death or migration is allowed we need extra conditions to ensure that there is a uniform bound on the total rate at which these occur within a given cube $C_v$.

In some examples there is a hard-core constraint preventing particles from appearing at a distance less than $\varepsilon$ from each other, for some fixed $\varepsilon$. In this case there is a uniform bound on the number of particles that can coexist in $C_v$, and provided the migration and death rates for each particle are bounded the total rate at which events happen in $C_v$ will be bounded.

On the other hand, if the total number of potential particles coexisting in a cube $C_v$ is unbounded, for (2.4) to hold we need some stronger condition on the death and migration rates per particle than mere boundedness, since the jump rate at site $v$ includes the sum of death and migration rates over *all* particles in $C_v$. This could happen, for example, if the presence of a large number of particles nearby slows down the rate of updating of an individual particle in some suitable manner.

Examples of the continuum particle system setup just described are numerous. Markov processes in $\mathcal{S}$ are discussed in [26; 27; 33; 52], and our results add to these works. Applications of these processes include the simulation of point processes in statistical literature (see e.g. [3; 23; 43]). In the case of deposition models, the 'state' of a particle in the set $\mathcal{T}$ could represent, for example, its shape, and/or its height above the surface.

The following sections contain discussions of some examples. Our results (Theorems 2.1, 3.1 and 3.2) provide a systematic framework to approach existence, spatial LLNs and spatial CLTs for these examples, although it should be pointed out that results along the lines of the LLN (Theorem 3.1) and our first CLT (Theorem 3.2, at least for one-dimensional distributions) can also be achieved for some of these examples by ad-hoc modifications to the arguments in Penrose and Yukich [50]; indeed, this was shown by J.E. Yukich (personal communication 2001), who suggested many of these examples to the author.

### *4.2. Random Sequential Adsorption and variants*

In the simplest version of Random Sequential Adsorption (RSA), particles do not carry marks and the only events are immigrations; there is no migration, death, or synchronous updating. An incoming particle is accepted with probability 1 if no existing particle lies within unit distance of that particle, and 0 otherwise. In effect, each particle lies at the centre of a hard sphere of unit diameter, and is accepted if its associated sphere does not overlap any previously accepted hard sphere.

The model just described fits our framework, provided we take $R \geq 1$, and the existence result (Theorem 2.1) shows this process is well-defined in $\mathbb{R}^d$. Our results give LLNs and CLTs for additive functionals of the form $S_\phi^{\tilde{A}}$ on $\mathcal{S}$, applied to the process $(\mathcal{X}_t^{\mathbf{A}})$, for which particles arrive only over the bounded region $\mathbf{A} \subset \mathbb{R}^d$. Note that in the basic model, $\phi$ is effectively a function on $\mathcal{S}_R$ rather than on $\mathcal{T} \times \mathcal{S}_R$, because points are unmarked. Examples include the following:

- If we take $\phi_1(\mathcal{X}) = 1$ for all $\mathcal{X} \in \mathcal{S}_R$, then $S_{\phi_1}^{\mathbf{A}}(\mathcal{X})$ (for $\mathcal{X} \in \mathcal{S}$) simply counts the number of accepted points in $\mathbf{A}$.
- Let $R_1$ be a parameter with $1 \leq R_1 \leq R$, and let $\mathrm{card}(\mathcal{X})$ denote the number of elements of $\mathcal{X}$. If (for $\mathcal{X} \in \mathcal{S}_R$) we take $\phi_2(\mathcal{X}) = (1/2)\mathrm{card}(\mathcal{X} \cap$ $B_{R_1}$), then $S^{\mathbf{A}}_{\phi_2}(\mathcal{X})$ (defined for $\mathcal{X} \in \mathcal{S}$) counts the number of pairs of points within distance $R_1$ of each other with both points in $\mathbf{A}$, plus half the number of such pairs with one point in $\mathbf{A}$.

Generalizations of RSA include cooperative sequential adsorption (where the probability of acceptance is some function of the local configuration of points), RSA with the identical hard spheres of the basic model replaced by independent identically distributed randomly shaped hard objects (so particles carry marks, representing their shape), and RSA with desorption where particles leave (die) at a non-zero rate. RSA and variants are of widespread interest in applications, as discussed in [50], and references therein.

Using Theorems 3.1 and 3.2, we may obtain LLNs and CLTs for $S^{\tilde{A}_n}_{\phi_i}(\mathcal{X}^{\tilde{A}_n,\nu}_t)$ with $i = 1$ or $i = 2$. In the case of $\phi_1$, these have previously been given in [50]. Theorem 3.3 can be used to extend the CLT to a CLT in $D[0, \infty)$ for $(S^{\tilde{A}_n}_{\phi_i}(\mathcal{X}^{\tilde{A}_n,\nu}_t), t \geq 0)$.

### 4.3. Continuum ballistic deposition

Ballistic deposition (BD) models in the continuum have the following features. Particles are represented, typically, by unit diameter hard balls in $d+1$ dimensions (usually $d = 1$ or $d = 2$). They fall vertically from above, sequentially at random, towards the adsorption surface $\mathbb{R}^d$ (strictly speaking, $\mathbb{R}^d \times \{0\}$). If a particle reaches the surface unhindered, it is attached to it. Details of what happens if a particle strikes another particle before reaching the surface depend on the model.

In *monolayer* versions of BD, all accepted particles lie on the surface of the substrate, as in the basic RSA model, so necessarily some particles are rejected. In the version proposed for $d = 1$ by Solomon [59] and for $d = 2$ by Jullien and Meakin [37], on striking a previously deposited particle, the new particle rolls, following the path of steepest descent until it reaches a stable position. If this position touches the adsorption surface the particle is fixed there; otherwise the particle is rejected. It is assumed that the rolling occurs instantaneously.

All accepted particles lie on the substrate, and so can be represented by points in $\mathbb{R}^d$. The position of an accepted particle is a translate (or *displacement*) of the location in $\mathbb{R}^d$ above which it originally comes in. There is a uniform bound $R_2$ on the size of the possible displacements (this is clear when $d = 1$, and a proof for $d = 2$ is given in [46]; see also [13]). Since the displacement and the decision on whether to accept an incoming particle at $x$ are determined by the existing configuration within distance $R_2 + 1$ of $x$, this model fits into our framework. There are many other monolayer BD models satisfying this condition of uniformly bounded displacements. Again, the functionals $S^{\mathbf{A}}_{\phi_1}(\mathcal{X})$ and $S^{\mathbf{A}}_{\phi_2}(\mathcal{X})$ are of interest; for $R_1 = 1$ the latter counts the number of touching pairs of particles in $\mathbf{A}$.

Senger et al. [56] describe many experimental results. Choi et al. [13] give both experimental and analytical results. Penrose [46] gives rigorous results on the infinite input version of this model.

In *multilayer* BD, a particle (ball) may attach itself to previously adsorbed particles instead of to the substrate. In the simplest form of continuum multilayer BD, each particle falls vertically towards the substrate and as soon as it encounters either the substrate or another particle, it sticks (and remains in that place forever).

In an alternative model of multilayer BD, when an incoming ball strikes a deposited ball, it does not stick but rolls, following the path of steepest descent, until it reaches a stable position. At this point it stops.

In either of these versions of multilayer BD, each particle is accepted, and has a vertical displacement (or *height*) relative to the position it would occupy if it were to fall to the substrate unhindered by other particles. If we take $\mathcal{T} = [0, \infty)$, this fits into our framework, and with the mark of a point representing the height above the substrate at which it is attached. In the case of multilayer BD with rolling, to apply our results we need a uniform bound on the horizontal distance a particle can roll.

### *4.4. Deposition models with rearrangement*

The RSA and BD models can be made arguably more realistic if an incoming particle is allowed to cause some (instantaneous, finite range) displacement of existing particles. If we modify the RSA and monolayer BD models in this way, our general existence results, and limit theorems for functionals such as $S^{\mathbf{A}}_{\phi_1}(\mathcal{X}^{\mathbf{A}}_t)$ and $S^{\mathbf{A}}_{\phi_2}(\mathcal{X}^{\mathbf{A}}_t)$ remain applicable, thereby adding to results in [50] which did not allow for synchronous updating, and to known experimental results for monolayer ballistic deposition with restructuring; see e.g. [37].

In the case of multilayer BD, such displacements can be viewed as 'avalanches'. Again, our general results remain applicable, provided we assume displacement has a finite range. Such an assumption is consistent with the assumption of high friction or bonding forces (see Onoda and Liniger [44]).

Models of this type are related to the notion of *random loose sphere packing* (RLSP). This has long been of interest in Materials Science, since the classic experimental works of Bernal, Mason and Scott [9; 10; 55]. See Cumberland and Crawford [14] for a survey and over 300 references. While there is no one commonly agreed upon mathematical model for RLSP, our multilayer BD model, with finite-range rolling and displacement allowed, has many of the features of the physical RLSP experiments.

Among the additive set functions $S^{\mathbf{A}}_\phi(\mathcal{X})$ of interest are those generated by the following choices of $\phi$:

1. *Total number of balls up to a specified height.* Let $R_3 > 0$, and for $\mathcal{X} \in \mathcal{S}_R$ let $\phi_3(\theta, \mathcal{X})$ be equal to 1 if $\theta \leq R_3$, and 0 otherwise. Then $S^{\mathbf{A}}_{\phi_3}(\mathcal{X})$ (defined for $\mathcal{X} \in \mathcal{S}$) counts the number of points of $\mathcal{X}$ in $\mathbf{A}$ which have height at most $R_3$ above the substrate. If we take $R_3 = 1/2$, then $S^{\mathbf{A}}_{\phi_3}(\mathcal{X})$ gives the number of balls touching the substrate, since the balls have radius $1/2$.

For large $R_3$, $S^{\mathbf{A}}_{\phi_3}(\mathcal{X})$ gives an estimate for the density of packed particles. Our results are applicable to give a LLN and CLT for $S^{\tilde{A}_n}_{\phi_3}(\mathcal{X}^{\tilde{A}_n,\nu}_t)$ when $(A_n, n \geq 1)$ satisfies (3.4) and (3.5).

2. *Total number of contacts ('coordination number')*. Let $\phi_4(\theta, \mathcal{X})$ denote half the number of points $(\mathbf{x}, \theta') \in \mathcal{X}$ such that $\|(\mathbf{x}, \theta') - (\mathbf{0}, \theta)\| = 1$, where $\|\cdot\|$ denotes $(d+1)$-dimensional Euclidean distance. Then, modulo boundary effects, $S^{\mathbf{A}}_\phi(\mathcal{X})$ denotes the total number of contacts between particles in $\mathcal{X}$. Our general results give us

$$\lim_{n\to\infty} S^{\tilde{A}_n}_{\phi_4}(\mathcal{X}^{\tilde{A}_n,\nu}_t)/|\tilde{A}_n| = \kappa(t) \quad a.s. \tag{4.1}$$

with a corresponding CLT and functional CLT.

3. *Surface growth.* In the multilayer BD process, possibly with rearrangement, the 'active zone' or 'interface' of the resulting agglomeration may be defined as the set of exposed particles, with a particle termed 'exposed' if there is a positive chance that the next incoming particle could strike it. Let $\phi_5(\theta, \mathcal{X})$ be the height $\theta$ if a particle at $(\mathbf{0}, \theta)$ is exposed, and let $\phi_5(\theta, \mathcal{X})$ be zero otherwise. Then $S^{\mathbf{A}}_{\phi_5}(\mathcal{X})$ is the total height of exposed particles over $\mathbf{A}$. We can apply Theorems 3.1, 3.2 and 3.3 to get a LLN, CLT and functional CLT for $S^{\tilde{A}_n}_{\phi_5}(\mathcal{X}^{\tilde{A}_n,\nu}_t)$. This adds to results in [51], where updating in the BD model is not assumed and the functional CLT is not considered.

Estimates for coordination numbers are found in Bernal and Mason [10], Gervois et al [29] and Gotoh and Finney [30] (p. 202) and references 148-153 of [14]. Jodrey and Tory [36] describe computer simulations of a loose sphere packing model. Non-rigorous work of Gervois et al [29] (p. 2124) suggests a coordination number of $6.1 \pm 0.3$ for a particular version of RLSP model. Gotoh and Finney [30] suggest a coordination number of 6.0.

## *4.5. Off-Lattice interacting particles*

We indicate some continuum analogs of classic interacting particle systems on the lattice.

a. ('voter model' I) $\mathcal{T}$ represents a set of two 'colours'. There is only immigration, with no death or migration. When a point arrives, it picks a point at random (within a distance $R$ of it) and chooses its colour with probability $p$, $0 \leq p \leq 1$. If no existing point lies within distance $R$ of the new point, its colour is selected at random. If $\phi$ is the indicator of a particle having a given colour, then $S^{\mathbf{A}}_\phi(\mathcal{X}_t)$ is the number of particles in $\mathbf{A}$ of that colour.

b. ('Voter Model' II) When a point arrives, its adopts the colour of the nearest neighbour within distance $R$. If no such neighbour exists, it acquires its colour randomly. As the number of particles increases, the scale at which interactions takes place decreases.

c. ('Exclusion process') Suppose particles do not arrive or depart, but jump in the state space. Suppose $\beta(d\mathbf{x})$ is a probability measure on $B_R$, and $\lambda > 0$, $\varepsilon > 0$ are constants. Suppose a particle at $\mathbf{x} \in \mathbb{R}^d$ attempts to jump at rate $\lambda$, and when an attempt is made to jump, the destination of the attempted jump lies in $\mathbf{x}+\Gamma$ with probability $\beta(\Gamma)$. The jump is successful if no previous particle lies within distance $\varepsilon$ of the destination. If the jump is unsuccessful the particle remains at $\mathbf{x}$.

d. ('Zero range process') A particle at $x \in \mathbb{R}^d$ jumps at rate $\lambda_n$ if there are $n$ other particles distant at most $\varepsilon$ from $x$. When it jumps, its destination lies in $x+\Gamma$ with probability $\beta(\Gamma)$. Here $\varepsilon$ and $\beta$ are as in the preceding example. If the parameters $\lambda_n$ satisfy $\sup_n(n\lambda_n) < \infty$ then the bounded jump rate condition (2.4) holds.

## 5. Further examples and discussion

### 5.1. Discrete examples

Not all applications of the general results of Sections 2 and 3 can be expressed in terms of marked point process in the continuum. Here we briefly discuss some cases where $X$ is a discrete space.

When $X$ is finite, our limit theorems remain applicable. Indeed, if $X$ is finite then any choice of $H$ will satisfy (3.3) since $H$ can take only finitely many values on a finite set. However, we concentrate here on examples where $X$ is the infinite set $\{0, 1, 2, \ldots\}$.

A model of *multilayer lattice BD* is be defined as follows. Given a vertex set $V$ satisfying (2.2), assume particles arrive at each site $v \in V$ as a Poisson process of intensity $\lambda$. Each particle attaches itself to the system at a height which is 1 plus the maximum previous height of particles at sites in $\mathcal{N}_v$. The particle system recording the maximum particle height at each site is a $X^V$-valued process which jumps from state $\eta$ to state $\eta^v$ at rate $\lambda$ for each vertex $v \in V$, where

$$\eta^v(u) = \begin{cases} \max_{w \in \mathcal{N}_v} \eta(w) + 1 & \text{if } u = v \\ \eta(u) & \text{otherwise} \end{cases}$$

This type of model is of much interest in the physical sciences; see for example [5]. Mathematical studies [2; 57] have mostly been concerned with the one-dimensional case where $V$ is $\mathbb{Z}$ or a sub-interval thereof. Our Theorem 2.1 shows the process is well-defined on more general graphs. Moreover, it can be shown that for $V = \mathbb{Z}^d$, when $\xi_0 \equiv 0$, and the additive set function $S_H^A(\xi)$ is taken to be the sum $\sum_{v \in A}(\xi(v))^k$, for arbitrary $k$, then the moments condition (3.10) holds for all $\gamma$ and all $\tau$.

A modification of the BD model is to allow *surface relaxation* (see [5]). Again, $\xi_t(v)$ denotes the height of the highest particle at site $v$ at time $t$. But now an incoming particle at site $v$ seeks the site $w$ in $\mathcal{N}_v$ with the lowest height $\xi_t(w)$ (in the event of a tie, select $w$ at random from the minima) and locates itself there, adding 1 to the value of $\xi_t(w)$. Again, our results ensure existence of $\xi_t(w)$, and LLNs and CLTs for additive set functions of the form $\sum_{v \in A}(\xi(v))^k$.

Superficially related to these deposition models are 'sand-pile' type models, which have been much studied in recent years (see e.g. [17; 42]). However, they do not fit so easily into our framework, since there is typically no bound on the possible range of influence of a single incoming particle.

## 5.2. Open problems

Our general model has some limitations which restrict its scope; it would be of interest to remove some of these.

One limitation already mentioned is the uniform bound on jump rates. This rules out direct application to a number of spatial immigration, death and migration type models where there is no uniform bound on the total rates of migrations/deaths within a bounded region.

Another limitation already mentioned is the restriction to finite-range interactions, which prevents our results being directly applicable to sand-pile type models, or to cooperative sequential adsorption models where the probability of acceptance of a particle depends on the whole existing configuration, not just within a finite range.

In continuum spatial models, it is natural to try to allow for non-jumpwise motion. The simplest extension would be to allow for deterministic evolution of the system in between the jump-type events, with the deterministic evolution at one site affecting neighboring sites only through the jump rates. Such an extension could be sufficient to deal with models where the lifetime of a particle is not exponentially distributed; for a particular example of this see the particle system in [58]. The theory of piecewise deterministic Markov processes [15; 35] could be relevant here. In a more general extension, one might hope to allow piecewise deterministic motion of a particle. In this case, this motion is typically able to take the particle from one patch of space into a different patch of space (this feature does not arise in [58]), which is likely to make matters more complicated.

## 6. Proof of Theorem 2.1

For $v \in V$, and $x \in X^{\mathcal{N}_v}$, recall that $\alpha_v^*(x, \cdot)$ denotes the jump rate kernel of the Markov process on $X^{\mathcal{N}_z}$ with generator $G_v^*$, and recall the definition (2.4) of $c_{\max}$ as the supremum of all jump rates. As a first step, note that for $x \in X^V$, by (2.3) and (2.6) followed by (2.5) we have

$$G_v f(x) = \int_{X^{\mathcal{N}_v}} (f_x^v(y) - f_x^v(x|_{\mathcal{N}_v})) \alpha_v^*(x|_{\mathcal{N}_v}, dy)$$
$$= \int_{X^{\mathcal{N}_v}} (f(x|v|y) - f(x)) \alpha_v^*(x|_{\mathcal{N}_v}, dy). \qquad (6.1)$$

For finite $A \subset V$, let $\mathcal{C}(A)$ be the set of all bounded measurable functions $f : X^V \to \mathbb{R}$ which depend only on co-ordinates in $A$, i.e. which satisfy

$$f(x) = f(y), \ \forall x, y \in X^V \text{ such that } x(v) = y(v) \ \forall v \in A.$$

For finite $A \subset V$, let $\mathcal{N}_A := \cup_{v \in A} \mathcal{N}_v$. For $f \in \mathcal{C}$, let $\|f\| := \sup\{f(x) : x \in X^V\}$. We shall use the following result to establish uniqueness of the Markovian family of processes with generator $G$.

**Lemma 6.1.** *Suppose $A \subset V$ is finite and $f \in \mathcal{C}(A)$. Then*

$$\|G_v f\| \leq 2 c_{\max} \|f\|, \tag{6.2}$$

$$G_v f \equiv 0, \quad \forall v \in V \setminus \mathcal{N}_A, \tag{6.3}$$

*and*

$$G_v f \in \mathcal{C}(A \cup \mathcal{N}_v), \quad \forall v \in \mathcal{N}_A \tag{6.4}$$

*Proof.* Since the total measure of $\alpha_v^*(x|_{\mathcal{N}_v}, \cdot)$ is at most $c_{\max}$, (6.1) gives us (6.2). Also, if $v \in V \setminus \mathcal{N}_A$ and $y \in X^{\mathcal{N}_v}$, then the restriction of $x|v|y$ to $A$ is the same as the restriction of $x$ to $A$, so that $f(x|v|y) = f(x)$, and hence $G_v f(x) = 0$ by (6.1). This yields (6.3).

Now suppose $v \in \mathcal{N}_A$, and suppose $x, x' \in X^V$ with $x|_{A \cup \mathcal{N}_v} \equiv x'|_{A \cup \mathcal{N}_v}$. Then $x$ agrees with $x'$ on $\mathcal{N}_v$, so that the measures $\alpha_v^*(x|_{\mathcal{N}_v}, \cdot)$ and $\alpha_v^*(x'|_{\mathcal{N}_v}, \cdot)$ are identical. Moreover, $x$ agrees with $x'$ on $A$ so that $f(x') = f(x)$ and moreover for all $y \in X^{\mathcal{N}_v}$ we have that $f(x|v|y) = f(x'|v|y)$. Hence, by (6.1) we have

$$G_v f(x) = \int_{X^{\mathcal{N}_v}} (f(x'|v|y) - f(x')) \alpha_v^*(x'|_{\mathcal{N}_v}, dy) = G_v f(x').$$

Thus, $G_v f(x)$ depends on $x$ only through $x(w), w \in A \cup \mathcal{N}_v$, i.e. (6.4) holds. ∎

On a suitable probability space, let $(\mathcal{P}_v, v \in V)$ be a family of independent and identically distributed homogeneous Poisson point processes of intensity $c_{\max}$ in $[0, \infty) \times [0, 1]$. Label the points of $\mathcal{P}_v$ as $((T_i(v), U_i(v)), i = 1, 2, 3, \ldots)$, with $T_1(v) < T_2(v) < T_3(v) < \cdots$ (with probability 1 the $T_i(v)$ are distinct for each $v$).

For $A \subseteq V$ and $v \in V$ let $\mathcal{N}_A^+$ denote the 2-neighbourhood of $A$ consisting of all vertices a graph distance at most 2 from some element of $A$ in the graph $\mathcal{V}$, and let $\mathcal{N}_v^+$ be the two-neighbourhood of $v$, i.e. let

$$\mathcal{N}_A^+ := \mathcal{N}_{\mathcal{N}_A}; \quad \mathcal{N}_v^+ := \mathcal{N}_{\{v\}}^+ = \mathcal{N}_{\mathcal{N}_v}. \tag{6.5}$$

By (2.2), the number of elements of $\mathcal{N}_v^+$ is at most $1 + D + D(D-1) = 1 + D^2$. Define the point set $\mathcal{P} \subset V \times [0, \infty)$ by

$$\mathcal{P} := \cup_{v \in V} \left( \{(v, 0)\} \cup \{(v, T_i(v)) : i \geq 1\} \right). \tag{6.6}$$

Make $\mathcal{P}$ into the vertex set of an (infinite) oriented graph $\mathcal{G}$ by putting in an oriented edge $(u, T) \to (v, T')$ whenever $v \in \mathcal{N}_u^+$ and $T < T'$.

For $w, v \in V$, let us say that *$w$ affects $v$ before time $t$* if there exists a (directed) path in the oriented graph that starts at $(w, 0)$ and ends at some Poisson point $(v, T)$ with $T \leq t$. Let $E_t(w, v)$ denote the event that $w$ affects $v$ before time $t$. Let $\mathrm{dist}(\cdot, \cdot)$ denote graph distance.

**Lemma 6.2.** *There exists a constant $\delta > 0$, such that for all $v \in V$, and all $n \in \mathbb{N}$,*

$$P\left[\bigcup_{w \in V: \text{dist}(w,v) \geq 2n} [E_{\delta n}(w, v)]\right] \leq 2^{-n}; \qquad (6.7)$$

$$P\left[\bigcup_{w \in V: \text{dist}(w,v) \geq 2n} [E_{\delta n}(v, w)]\right] \leq 2^{-n}. \qquad (6.8)$$

*Proof.* Let $W$ be exponentially distributed with mean $1/c_{\max}$, and let $S_n$ be the sum of $n$ independent copies of $W$. Let $\theta > 0$. Then by similar arguments to those used in the proof of Lemma 5.1 of [47], the probability on the left side of (6.7) is bounded by $(D^2)^n P[S_n \leq \delta n]$, and hence by $D^{2n} e^{\theta \delta n} (E[e^{-\theta W}])^n$.

Choose $\theta$ so that $E[e^{-\theta W}] < 1/(4D^2)$, and then choose $\delta > 0$ so that $e^{\theta \delta} < 2$. Then the preceding bound is at most $(1/2)^n$, completing the proof of (6.7). The proof of (6.8) is similar. ∎

**Corollary 6.1.** *Let $v \in V$ and $t > 0$. Then* (a) *with probability 1, $v$ affects only finitely many $w \in V$ before time $t$, and* (b) *with probability 1, only finitely many $w \in V$ affect $v$ before time $t$.*

*Proof.* Let $F_n$ be the event that there exists $w \in V$ with $\text{dist}(w, v) \geq 2n$ such that $v$ affects $w$ before time $t$. Then $P[F_n] \to 0$ as $n \to \infty$ by (6.8), and this implies part (a). The proof of (b) using (6.7) is similar. ∎

We partition the vertex set $\mathcal{P}$ of $\mathcal{G}$ into 'generations' by putting each point $(v, 0), v \in V$ in generation 0 and saying a point $(v, T)$, $T \in (0, \infty)$ of $\mathcal{P}$ is in the $k$th generation if the *longest* directed path from any point $(u, 0)$ to $(v, T)$ is of length $k$ (with length measured by graph distance). By Corollary 6.1, this is indeed a partition with probability 1.

Recall that $\alpha_v^*(x)$ denotes the total measure of the jump rate kernel, and that $\alpha_v^*(x) \leq c_{\max}$. Let $\delta_x(\cdot)$ denote the unit point mass at $x$ and define the probability kernel

$$\mu_v(x, \cdot) = \frac{\alpha_v^*(x, \cdot) + (c_{\max} - \alpha_v^*(x))\delta_x(\cdot)}{c_{\max}}. \qquad (6.9)$$

For each $v \in V$, there is a measurable function $\psi_v : X^{\mathcal{N}_v} \times [0, 1] \mapsto X^{\mathcal{N}_v}$ with the property that if $\vartheta$ is uniformly distributed over $[0, 1]$, then $\psi_v(x, \vartheta)$ has distribution $\mu_v(x, \cdot)$ for each $x \in X^{\mathcal{N}_v}$. The existence of such a $\psi_v$ follows by Kallenberg ([38], page 56, lemma 3.22), since $\mu_v(x, B)$ is a probability kernel on $X^{\mathcal{N}_v}$, which is a Polish space, and hence a Borel space (see [38] theorem A1.2, or [24], p. 419, Proposition 20).

For $v \in V$, set $T_0^*(v) = 0$ and list the arrival times $T_j(u), u \in \mathcal{N}_v$ in increasing order as $T_1^*(v), T_2^*(v), T_3^*(v), \ldots$.

Given $x \in X^V$, construct the process $(\xi_t^x, t \geq 0)$ as follows. First, for all $v \in V$ assume that $\xi_t^x(v)$ is a right-continuous function of $t$ that is constant in between arrival times $T_i^*(v), i \geq 1$. In other words, for all $v \in V$ and $t > 0$, if $T_{i-1}^*(v) \leq t < T_i^*(v)$ then set $\xi_t^x(v) = \xi_{T_{i-1}^*}^x(v)$. Set $\xi_0^x(v) = x(v)$ for all $v \in V$. Suppose inductively that $\xi_T^x(v)$ is defined for each $(v, T)$ in generation $k$ of $\mathcal{P}$.

Suppose arrival $(v, T_j(v))$ is in generation $k + 1$. Set $T = T_j(v)$. Define the new state of the process in $\mathcal{N}_v$ at time $T$ by

$$\xi_T^x|_{\mathcal{N}_v} = \psi_v(\xi_{T-}^x|_{\mathcal{N}_v}, U_j(v)). \tag{6.10}$$

This is well-defined because for each $w \in \mathcal{N}_v$, there exists $i = i(w, T)$ such that $T = T_i^*(w)$, and $(w, T_{i-1}^*(w))$ is in generation $k$ at the latest, so that $\xi_t^x(w)$ is already defined for $t$ in the interval $[T_{i-1}^*(w), T)$ and hence $\xi_{T-}^x(w)$ is well-defined.

**Lemma 6.3.** *There exists a filtration $(\mathcal{G}_t, t \geq 0)$ with respect to which the family of processes $((\xi_t^x, t \geq 0), x \in X^V)$ described above is a Markovian family of processes in state-space $X^V$. The associated family of transition distributions induces a semigroup of operators $(P_t, t \geq 0)$ with generator $G$ given on $\mathcal{C}$ by (2.7); that is, with this $P_t$ and $G$, (2.1) holds with uniform convergence for all $f \in \mathcal{C}$.*

*Proof.* Let $\mathcal{G}_t$ denote the $\sigma$-field generated by all Poisson arrivals before time $t$ at all sites, i.e. let $\mathcal{G}_t$ denote the smallest $\sigma$-field with respect to which each of the random variables of the form $\mathcal{P}_v(B)$, with $v \in V$ and $B$ a Borel subset of $[0, t] \times [0, 1]$, is measurable. Here $\mathcal{P}_v(A)$ denotes the number of points in $A$ of the Poisson process $\mathcal{P}_v$.

By the construction and standard properties of the Poisson process, the processes $(\xi_\cdot^x, x \in X^V)$ form a Markov family of processes, adapted to $(\mathcal{G}_t, t \geq 0)$. Sample paths are right-continuous with left limits because we use the product topology and there are no accumulations of Poisson arrivals at any site $v \in V$.

Let $(P_t, t \geq 0)$ be the semigroup associated with the Markov family $(\xi_\cdot^x, x \in X^V)$. Let $f \in \mathcal{C}$ and let $x \in X^V$. Then

$$P_t f(x) - f(x) = E[f(\xi_t^x) - f(x)].$$

Let $x \in X^V$. Choose finite $A \subset V$ such that $f \in \mathcal{C}(A)$. Then $f(\xi_t^x) = f(x)$ unless $\min\{T_1(v), v \in \mathcal{N}_A\} \leq t$.

Let $F_t^{(1)}$ be the event that that $T_1(v) \leq t$ for two or more $v \in \mathcal{N}_A^+$, and let $F_t^{(2)}$ be the event that $T_2(v) \leq t$ for some $v \in \mathcal{N}_A$. Let $F_t := F_t^{(1)} \cup F_t^{(2)}$. Also, for $v \in \mathcal{N}_A$, let $F_{t,v}^{(3)}$ be the event that $T_1(v) \leq t < T_2(v)$ and $T_1(w) > t$ for $w \in \mathcal{N}_A^+ \setminus \{v\}$. Then

$$E[f(\xi_t^x) - f(x) | F_{t,v}^{(3)}] = \int_{X^{\mathcal{N}_v}} (f(x|v|y) - f(x)) \mu_v(x|_{\mathcal{N}_v}, dy)$$

$$= c_{\max}^{-1} \int_{X^{\mathcal{N}_v}} (f(x|v|y) - f(x)) \alpha_v^*(x|_{\mathcal{N}_v}, dy),$$

so that by (6.1),

$$E[f(\xi_t^x) - f(x) | F_{t,v}^{(3)}] = c_{\max}^{-1} G_v f(x).$$

Hence,

$$\begin{aligned}
&Ef(\xi_t^x) - f(x) \\
&= P[F_t]E[f(\xi_t^x) - f(x)|F_t] + \sum_{v \in \mathcal{N}_A} P[F_{t,v}^{(3)}]E[f(\xi_t^x) - f(x)|F_t^{(3)}] \\
&= P[F_t]E[f(\xi_t^x) - f(x)|F_t] + \sum_{v \in \mathcal{N}_A} P[F_{t,v}^{(3)}]c_{\max}^{-1}G_v f(x). \quad (6.11)
\end{aligned}$$

Using (6.11), then (6.3) and then (6.2), we obtain

$$\left|t^{-1}(Ef(\xi_t^x) - f(x)) - Gf(x)\right|$$

$$= \left|t^{-1}P[F_t]E[f(\xi_t^x) - f(x)|F_t] + \sum_{v \in \mathcal{N}_A}(t^{-1}c_{\max}^{-1}P[F_{t,v}^{(3)}] - 1)G_v f(x)\right|$$

$$\leq 2\|f\|\left(t^{-1}P[F_t] + \sum_{v \in \mathcal{N}_A}\left(t^{-1}c_{\max}^{-1}P[F_{t,v}^{(3)}] - 1\right)\right). \quad (6.12)$$

As $t \downarrow 0$ we have that $P[F_t] = o(t)$ and $P[F_{t,v}^{(3)}] = c_{\max}t + o(t)$; moreover $P[F_{t,v}^{(3)}]$ is the same for all $v \in \mathcal{N}_A$. Hence the bound in (6.12) is independent of $x$ and tends to 0 as $t \downarrow 0$. Thus the generator of the semigroup $(P_t, t \geq 0)$ is $G$ on functions $f \in \mathcal{C}$. ∎

The next step is to show that the transition semigroup $(P_t, t \geq 0)$ satisfies (2.8). For finite $A \subset V$, to get a process with generator $G^A$ rather than $G$, we need to 'switch off' all jumps taking places at sites outside $A$. To do this, set

$$\mathcal{P}_A := \mathcal{P} \cap (A \times [0, \infty)),$$

and partition $\mathcal{P}_A$ into generations in precisely the same manner as for $\mathcal{P}$. Then using the updating rule (6.10) but now restricting attention to arrival times $T_j(v), v \in A$, we obtain a pure jump process $(\xi_t^{A,x}, t \geq 0)$ with generator $G_A$. Note that we use the *same* family of Poisson processes $(\mathcal{P}_v, v \in V)$ in constructing $(\xi_t^{A,x}, t \geq 0)$ and $(\xi_t^x, t \geq 0)$, regardless of $A$.

For each $(v, t) \in V \times [0, \infty)$, let the 'cluster' $C_{v,t}$ be the set of $w \in V$ such that $w$ affects $z$ before time $t$ for some $z \in \mathcal{N}_v^+$ (recall $\mathcal{N}_v^+$ is the 2-neighbourhood of $v$). By Corollary 6.1 this cluster is almost surely finite. Its significance lies in the following result, which will be used again in proving limit theorems later on.

**Lemma 6.4.** *Suppose* $(v, t) \in V \times [0, \infty)$, *and* $x \in X^V$, *and* $A \subset V$ *is finite. If* $C_{v,t} \subseteq A$, *then* $\xi_t^{A,x}(w) = \xi_t^x(w)$ *for all* $w \in \mathcal{N}_v$.

*Proof.* All differences between $(\xi_s^{A,x}, s \geq 0)$ and $(\xi_s^x, s \geq 0)$ are due to Poisson arrivals in $(V \setminus A) \times [0, \infty) \times [0, 1]$ which contribute to changes in $(\xi_s^x, s \geq 0)$ but not to changes in $(\xi_s^{A,x}, s \geq 0)$. However, by the very definition of $C_{v,t}$, since $C_{v,t} \subseteq A$ the influence of these Poisson arrivals does not propagate to any vertex in $\mathcal{N}_v^+$ by time $t$, and therefore $\xi_t^x|_{\mathcal{N}_v} \equiv \xi_t^{A,x}|_{\mathcal{N}_v}$. ∎

**Lemma 6.5.** *The transition semigroup $(P_t, t \geq 0)$ of the Markovian family of processes $((\xi_t^x, t \geq 0), x \in X^V)$ defined above satisfies (2.8).*

*Proof.* Suppose $f \in \mathcal{C}$ and choose finite $A \subset V$ such that $f \in \mathcal{C}(A)$. Suppose $(W_m, m \geq 1)$ is a sequence of finite subsets of $V$ satisfying $\liminf_{m \to \infty} W_m = V$. Then for $v \in A$, by Lemma 6.4 we have $\xi_t^{W_m,x}(v) = \xi_t^x(v)$ for large enough $m$, almost surely. Hence, $f(\xi_t^{W_m,x}) = f(\xi_t^x)$ for large enough $m$, almost surely. Hence by the dominated convergence theorem,

$$P_t f(x) = E[f(\xi_t^x)] = \lim_{m \to \infty} E[f(\xi_t^{W_m,x})]$$
$$\lim_{m \to \infty} P_t^{W_m}(x). \quad \blacksquare$$

*Proof of Theorem 2.1.* In view of Lemmas 6.3 and 6.5, it remains only to show that the transition semigroup $(P_t, t \geq 0)$ is the only one with generator $G$ on functions in $\mathcal{C}$. We do this directly, without appealing to generalities on resolvents, in part because the state space is not locally compact in general. Let $f \in \mathcal{C}$, and choose a finite set $A \subset V$ such that $f \in \mathcal{C}(A)$. Write $\mathcal{N}_A^+$ for $\cup_{v \in A} \mathcal{N}_v^+$. By Lemma 6.1,

$$Gf = \sum_{v \in \mathcal{N}_A} G_v f \in \mathcal{C}. \tag{6.13}$$

For $t \geq 0$ and $h > 0$, we have

$$\|h^{-1}(P_{t+h}f - P_t f) - P_t Gf\| = \left\|P_t\left(\left(\frac{P_h f - f}{h}\right) - Gf\right)\right\| \to 0 \quad \text{as } h \downarrow 0,$$

because by assumption $G$ is the generator, defined on $f \in \mathcal{C}$ by assumption so that (2.1) holds with uniform convergence, and $P_t$ is a contraction. Similarly, for $t > 0$ and $0 < h < t$,

$$\|h^{-1}(P_t f - P_{t-h}f) - P_t Gf\| = \left\|P_{t-h}\left(\left(\frac{P_h f - f}{h}\right) - Gf\right)\right\| \to 0 \quad \text{as } h \downarrow 0.$$

Hence for any $x \in X^v$, we have $\frac{d}{dt}P_t f(x) = P_t Gf(x)$. Since $Gf \in \mathcal{C}$ by (6.13), for integer $n > 0$ we can repeat the preceding argument $n$ times to deduce that

$$\frac{d^n}{dt^n}P_t f(x) = P_t G^n f(x).$$

Hence by Taylor's theorem with Lagrange remainder,

$$P_t f(x) = \left(\sum_{j=0}^n \frac{G^j f(x) t^j}{j!}\right) + \frac{P_u G^{n+1} f(x) t^{n+1}}{(n+1)!}, \quad \text{some } u \in (0, t). \tag{6.14}$$

By twice applying Lemma 6.1, we have that for $v \in \mathcal{N}_A$ and $w \in V$, if $w \notin \mathcal{N}_{A \cup \mathcal{N}_v}$ then $G_w G_v f \equiv 0$, whilst in any case $G_w G_v f \in \mathcal{C}(A \cup \mathcal{N}_v \cup \mathcal{N}_w)$. Hence, if $w \notin \mathcal{N}^+_{A \cup \{v\}}$ (as defined in (6.5)) then $G_w G_v f \equiv 0$, whilst in any case $G_w G_v f \in \mathcal{C}(A \cup \mathcal{N}_{\{v,w\}})$. Repeating the same argument, an induction on $n$ shows that $G_{v_n} G_{v_{n-1}} \cdots G_{v_1} f \equiv 0$ except when each $v_i$ lies in $\mathcal{N}^+_{A \cup \{v_1, \ldots, v_{i-1}\}}$. Hence, iterating (6.13) yields

$$G^n f(x) = \sum_{v_1 \in \mathcal{N}^+_A} \sum_{v_2 \in \mathcal{N}^+_{A \cup v_1}} \cdots \sum_{v_n \in A \cup \mathcal{N}^+_{A \cup \{v_1, \ldots, v_{n-1}\}}} G_{v_n} G_{v_{n-1}} \cdots G_{v_1} f(x). \quad (6.15)$$

Recall the definition (2.2) of $D$, and let $|\cdot|$ denote cardinality. The number of $n$-tuples $(v_1 \ldots, v_n) \in V^n$ such that each $v_i$ lies in $\mathcal{N}^+_{A \cup \{v_1, \ldots, v_{i-1}\}}$ is bounded by the expression

$$|A|(|A| + 1)(|A| + 2) \cdots (|A| + n - 1)(1 + D^2)^n$$
$$\leq 2^n \max(|A|, 1) \max(|A|, 2) \cdots \max(|A|, n)(1 + D^2)^n$$
$$\leq n! |A|^{|A|} 2^n (1 + D^2)^n.$$

Moreover, by (6.2), for all $v_1, v_2, \ldots, v_n$ we have $\|G_{v_n} \cdots G_{v_1} f\| \leq (2c_{\max})^n \|f\|$. Thus, setting $K_0 := 4c_{\max}(1+D^2)$, by (6.15) and the fact that $P_u$ is a contraction we have that

$$\|P_u G^n f\| \leq \|G^n f\| \leq n! |A|^{|A|} K_0^n \|f\|, \quad u \geq 0.$$

Thus for $t < 1/K_0$, the remainder term in (6.14) tends to zero as $n \to \infty$, so

$$P_t f(x) = \sum_{j=0}^{\infty} t^j G^j f(x)/j!, \quad t < 1/K_0.$$

Thus, for $f \in \mathcal{C}$ and $t < 1/K_0$, we have an explicit expression for $P_t(f)$ in terms of $G$. Hence, if $(\tilde{P}_t)_{t \geq 0}$ is another semigroup, associated with a Markov family of processes in $X^V$, with generator $G$ on $\mathcal{C}$, we must have $\tilde{P}_t(f) \equiv P_t(f)$ for all $t < 1/K_0$. By the Monotone Class Theorem [62], this identity extends to all bounded measurable functions on $X^V$, and by the semigroup property it then extends to all $t > 0$. This completes the proof of uniqueness. ∎

## 7. Proof of general limit theorems

To prove Theorems 3.1 and 3.2 we shall apply results from [47]. Recall that $V = \mathbb{Z}^d$ for these results and $\mathcal{B}$ denotes the collection of finite subsets of $\mathbb{Z}^d$. Given a probability distribution $\nu$ on $X$, let the process $(\xi_t^\nu, t \geq 0)$ and (for $A \in \mathcal{B}$) the process $(\xi_t^{A,\nu}, t \geq 0)$ be as described in Section 3.1. That is, $\xi_\cdot^\nu$ (respectively $\xi_\cdot^{A,\nu}$) has generator $G$ (respectively $G^A$) with initial distribution given by the product measure with all marginals equal to $\nu$.

We assume for our proofs that the process $(\xi_t^\nu, t \geq 0)$ is constructed in terms of the family $(\mathcal{P}_v, v \in \mathbb{Z}^d)$ of Poisson processes in $[0, \infty) \times (0, 1)$, as described in Section 6, in particular the description leading up to (6.10), with initial values $(\xi_0^\nu(v), v \in \mathbb{Z}^d)$ given by a family of independent $X$-valued variables $(\zeta_v, v \in \mathbb{Z}^d)$ with common distribution $\nu$, independent of the Poisson processes $\mathcal{P}_v, v \in \mathbb{Z}^d$.

Likewise we assume the process $(\xi_t^{A,\nu}, t \geq 0)$ is constructed using the same family of Poisson processes, as described just before the statement of Lemma 6.4, and using the same initial values $(\zeta_v, v \in \mathbb{Z}^d)$. Thus we couple the processes $\xi_\cdot^\nu$ and $\xi_\cdot^{A,\nu}$ so the initial values are identical $\xi_0^\nu = \xi_0^{A,\nu}$ and the same Poisson processes $(\mathcal{P}_v, v \in \mathbb{Z}^d)$ are used in both cases (but in the case of $\xi_t^A$ the effects of the Poisson processes $\mathcal{P}_v$ for $v \notin A$ are 'switched off').

Recall from Section 6 that for $(v, t) \in \mathbb{Z}^d \times [0, \infty)$ the 'cluster' $C_{v,t}$ is the set of $u \in \mathbb{Z}^d$ such that $u$ affects $v$ before time $t$ for some $v \in \mathcal{N}_v^+$, and is almost surely finite. Since with our coupling we have $\xi_0^\nu = \xi_0^{A,\nu}$, the following result is immediate from Lemma 6.4.

**Lemma 7.1.** *Suppose $(v, t) \in \mathbb{Z}^d \times [0, \infty)$ and suppose $A \subset \mathbb{Z}^d$ is finite. If $C_{v,t} \subseteq A$, then $\xi_t^{A,\nu}(w) = \xi_t^\nu(w)$ for all $w \in \mathcal{N}_v$.*

Let $t > 0$. Given $A \in \mathcal{B}$ and $v \in A$, $t > 0$, define $Y_v(A, t)$ by

$$Y_v(A, t) := \begin{cases} H(L_v(\xi_t^{A^o,\nu}|\mathcal{N}_v)) & \text{if } v \in A^o \\ 0 & \text{if } v \in A \setminus A^o \end{cases} \quad (7.1)$$

If $(\mathcal{N}_A)^o = A$, we have

$$S_H^A(\xi_\tau^{A,\nu}) = \sum_{v \in (\mathcal{N}_A)^o} H(L_v(\xi_\tau^{(\mathcal{N}_A)^o,\nu}|\mathcal{N}_v)) = \sum_{v \in \mathcal{N}_A} Y_v(\mathcal{N}_A, \tau).$$

In the terminology of Section 3 of [47], for any $t > 0$ the random variables $Y_z(A, t)$, defined for $z \in A, A \in \mathcal{B}$, form a *stationary $\mathcal{B}$-indexed summand* with respect to the i.i.d. family $((\mathcal{P}_z, \zeta_z), z \in \mathbb{Z}^d)$. To prove Theorems 3.1 and 3.2, we shall use Theorem 3.1 of [47]. We need to check the conditions of this result.

*Proof of Theorem 3.1.* If $(A_n, n \geq 1)$ is a $\mathcal{B}$-valued sequence with $\liminf A_n = \mathbb{Z}^d$, then with probability 1, $A_n$ includes the cluster $C_{\mathbf{0},\tau}$ for all large enough $n$, and hence by Lemma 7.1, $\xi_\tau^{A_n,\nu}(v) = \xi_\tau^\nu(v)$ for all but finitely many values of $n$, and all $v \in \mathcal{N}$. Hence, $Y_{\mathbf{0}}(A_n, \tau)$ tends to $Y_{\mathbf{0},\tau} := H(\xi_\tau^\nu|\mathcal{N})$ as $n \to \infty$, almost surely. Moreover, the variables $Y_{\mathbf{0}}(A_n, \tau)$ are uniformly integrable by the assumption that (3.3) holds for some $\gamma > 1$, and therefore the convergence of $Y_{\mathbf{0}}(A_n, \tau)$ to $Y_{\mathbf{0},\tau}$ extends to convergence in $L^1$.

Thus, all the conditions for the first part of Theorem 3.1 of [47] hold (note that $|A_n|/|\mathcal{N}_{A_n}| \to 1$ by (3.4)), and by (3.3) of [47] we obtain (3.6). ∎

We now work towards a proof of Theorem 3.2. Given $t > 0$, $v \in \mathbb{Z}^d$, and $A \in \mathcal{B}$, define the random variables

$$Y_{v,t} := H(L_v(\xi_t^\nu|\mathcal{N}_v)); \quad Y_{v,t}^{(A)} := H(L_v(\xi_t^{A_n,\nu}|\mathcal{N}_v)). \quad (7.2)$$

We now give some bounds on covariances. For $v \in \mathbb{Z}^d$ and $A \subset \mathbb{Z}^d$ with $A \neq \emptyset$, set $\mathrm{dist}(v, A) := \min\{|v - w| : w \in A\}$. Also, recall from Section 3.1 that $B_r$ is the closed Euclidean ball centred at the origin in $\mathbb{R}^d$, and for $r > 0$ and $v \in \mathbb{Z}^d$, write $B_r(v)$ for $v + B_r$.

**Lemma 7.2.** *Let $s > 0, t > 0$ and assume the moments condition (3.3) holds for $\tau \in \{s,t\}$, for some $\gamma > 2$. Then there is a constant $K$ such that for all $v \in \mathbb{Z}^d, w \in \mathbb{Z}^d$ and $A \in \mathcal{B}$ we have*

$$|\mathrm{Cov}(Y_{v,s}, Y_{w,t})| \leq K \exp(-|v-w|/K) \tag{7.3}$$

*and this also holds if $Y_{v,s}$ is replaced by $Y_{v,s}^{(A)}$ and/or $Y_{w,t}$ is replaced by $Y_{w,t}^{(A)}$, and also*

$$|\mathrm{Cov}(Y_{v,s}^{(A)} - Y_{v,s}, Y_{w,t}^{(A)})| \leq K \exp(-\mathrm{dist}(v, A^c)) \tag{7.4}$$

*and this still holds if $Y_{w,t}^{(A)}$ is replaced by $Y_{w,t}$.*

*Proof.* Let $v \in \mathbb{Z}^d, w \in \mathbb{Z}^d$. Let $r := (|w-v|/3)$. Let $\hat{Y}_{v,s} := Y_{v,s}\mathbf{1}_{\{C_{v,s} \subseteq B_r(v)\}}$. Then $\hat{Y}_{v,s}$ is determined by $((\mathcal{P}_z, \zeta_z), z \in B_r(v) \cap \mathbb{Z}^d)$, so that $\mathrm{Cov}(\hat{Y}_{v,s}, \hat{Y}_{w,t}) = 0$. Choose $\gamma > 2$ such that (3.3) holds for $\tau \in \{s,t\}$. By Hölder's inequality and Lemma 6.2, there exists a constant $K$ (independent of $v,w$) such that

$$\begin{aligned} E[(Y_{v,s} - \hat{Y}_{v,s})^2] &= E[Y_{v,s}^2 \mathbf{1}_{\{C_{v,s} \setminus B_r(v) \neq \emptyset\}}] \\ &\leq (E[Y_{v,s}^\gamma])^{2/\gamma} P[C_{v,s} \setminus B_r(v) \neq \emptyset]^{(\gamma-2)/\gamma} \\ &\leq K \exp(-r/K) \end{aligned}$$

and similarly, $E[(Y_{w,t} - \hat{Y}_{w,t})^2] \leq K\exp(-r/K)$ so that using the Cauchy-Schwarz inequality and the moments condition (3.3), we obtain for another constant (also denoted $K$) that

$$\begin{aligned} |\mathrm{Cov}(Y_{v,s}, Y_{w,t})| &\leq |\mathrm{Cov}(Y_{v,s} - \hat{Y}_{v,s}, Y_{w,t})| + |\mathrm{Cov}(\hat{Y}_{v,s}, Y_{w,t} - \hat{Y}_{w,t})| \\ &\leq K \exp(-r/K) \end{aligned}$$

so that (7.3) holds. The proof is identical with $Y_{v,s}$ replaced by $Y_{v,s}^{(A)}$ or $Y_{w,t}$ replaced by $Y_{w,t}^{(A)}$.

Turning to (7.4), put $r' = \mathrm{dist}(v, A^c) - 1$. Then $Y_{v,s} = Y_{v,s}^{(A)}$ unless $C_{v,t} \setminus B_{r'}(v) \neq \emptyset$, so by Hölder's inequality, (3.3) and Lemma 6.2,

$$\begin{aligned} E[(Y_{v,s}^{(A)} - Y_{v,s})^2] &\leq E[(Y_{v,s}^{(A)} - Y_{v,s})^\gamma]^{2/\gamma} P[C_{v,t} \setminus B_{r'}(v) \neq \emptyset]^{(\gamma-2)/\gamma} \\ &\leq K \exp(-r'/K), \end{aligned}$$

and combining this with the Cauchy-Schwarz inequality and using (3.3), we obtain (7.4). ∎

We now give a limit for covariances, which is part of the statement of Theorem 3.2.

**Lemma 7.3.** *Let $s \geq 0$ and $t \geq 0$. Suppose for some $\gamma > 2$ that (3.3) holds for $\tau \in \{s, t\}$. Suppose $(A_n)_{n \geq 1}$ is a $\mathcal{B}$-valued sequence satisfying (3.4) and (3.5). Define $\sigma(s, t)$ by (3.7). Then as $n \to \infty$,*

$$|A_n|^{-1} \mathrm{Cov}(S_H^{A_n}(\xi_s^{A_n, \nu}), S_H^{A_n}(\xi_t^{A_n, \nu})) \to \sigma(s, t). \tag{7.5}$$

*Proof.* Observe that $((Y_{v,s}, Y_{v,t}), v \in \mathbb{Z}^d)$ (defined by (7.2)) form a stationary random field. For $r > 0$, let $A_n^r$ be the set of sites $v \in \mathbb{Z}^d$ such that $B_r(v) \cap \mathbb{Z}^d \subseteq A_n$. Then for any constant $r > 0$,

$$|A_n|^{-1} \mathrm{Cov}(S_H^{A_n}(\xi_s^\nu), S_H^{A_n}(\xi_t^\nu)) = |A_n|^{-1} \sum_{v \in A_n} \sum_{w \in A_n} \mathrm{Cov}(Y_{v,s}, Y_{w,t})$$

$$= |A_n|^{-1} \sum_{v \in A_n^r} \left( \sigma(s, t) - \sum_{w \in \mathbb{Z}^d \setminus A_n} \mathrm{Cov}(Y_{v,s}, Y_{w,t}) \right)$$

$$+ |A_n|^{-1} \sum_{v \in A_n \setminus A_n^r} \sum_{w \in A_n} \mathrm{Cov}(Y_{v,s}, Y_{w,t}). \tag{7.6}$$

Every $v \in A_n \setminus A_n^r$ lies in $B_r(v')$ for some $v' \in \partial_{\mathrm{ext}}(A_n)$, so that $|A_n \setminus A_n^r| \leq |B_r \cap \mathbb{Z}^d| \times |\partial_{\mathrm{ext}}(A_n)|$, and so by (3.4), as $n \to \infty$ we have $|A_n \setminus A_n^r|/|A_n| \to 0$ so that $|A_n^r|/|A_n| \to 1$. Also, by (7.3), for some constant $K_1$,

$$|A_n|^{-1} \left| \sum_{v \in A_n \setminus A_n^r} \sum_{w \in A_n} \mathrm{Cov}(Y_{v,s}, Y_{w,t}) \right| \leq K_1 |A_n|^{-1} |A_n \setminus A_n^r| \to 0. \tag{7.7}$$

Also,

$$\limsup_{n \to \infty} |A_n|^{-1} \sum_{v \in A_n^r} \sum_{w \in \mathbb{Z}^d \setminus A_n} |\mathrm{Cov}(Y_{v,s}, Y_{w,t})| \leq \sum_{z \in \mathbb{Z}^d \setminus B_r} |\mathrm{Cov}(Y_{\mathbf{0},s}, Y_{z,t})| := h_1(r),$$

where $h_1(r) \to 0$ as $r \to \infty$ by (7.3). Combining this with (7.7) in (7.6), we obtain the limit

$$\lim_{n \to \infty} |A_n|^{-1} \mathrm{Cov}(S_H^{A_n}(\xi_s^\nu), S_H^{A_n}(\xi_t^\nu)) = \sigma(s, t). \tag{7.8}$$

Next, writing $Y_{v,t}^{(n)}$ for $Y_{v,t}^{(A_n)}$, we approximate the non-stationary random field $((Y_{v,s}^{(n)}, Y_{v,t}^{(n)}), v \in \mathbb{Z}^d)$ by the stationary random field $((Y_{v,s}, Y_{v,t}), v \in \mathbb{Z}^d)$. We write

$$\mathrm{Cov}(S_H^{A_n}(\xi_s^{A_n,\nu}), S_H^{A_n}(\xi_t^{A_n,\nu})) - \mathrm{Cov}(S_H^{A_n}(\xi_s^\nu), S_H^{A_n}(\xi_t^\nu))$$

$$= \sum_{v \in A_n} \sum_{w \in A_n} (\mathrm{Cov}(Y_{v,s}^{(n)}, Y_{w,t}^{(n)}) - \mathrm{Cov}(Y_{v,s}, Y_{w,t}))$$

$$= \left( \sum_{v \in A_n} \sum_{w \in A_n} \mathrm{Cov}(Y_{v,s}^{(n)} - Y_{v,s}, Y_{w,t}^{(n)}) \right)$$

$$+ \sum_{v \in A_n} \sum_{w \in A_n} \mathrm{Cov}(Y_{v,s}, Y_{w,t}^{(n)} - Y_{w,t}). \tag{7.9}$$

By Lemma 7.2, there is a constant $K_2$ such that

$$|\mathrm{Cov}(Y_{v,s}^{(n)} - Y_{v,s}, Y_{w,t}^{(n)})| \leq K_2 \exp(-\max(\mathrm{dist}(v, A_n^c), |v-w|)/K_2) \quad (7.10)$$

and from this we may deduce that for $v \in A_n^r$ we have

$$\sum_{w \in A_n} |\mathrm{Cov}(Y_{v,s}^{(n)} - Y_{v,s}, Y_{w,t}^{(n)})|$$
$$\leq \sum_{z \in \mathbb{Z}_d \setminus B_r} K_2 \exp(-|z|/K_2) + \sum_{z \in \mathbb{Z}^d \cap B_r} K_2 \exp(-r/K_2)$$
$$=: h_2(r)$$

where $h_2(r) \to 0$ as $r \to \infty$. Also, for any $r$ we have by (7.10) that for some constant $K_3$,

$$|A_n|^{-1} \sum_{v \in A_n \setminus A_n^r} \sum_{w \in \mathbb{Z}^d} |\mathrm{Cov}(Y_{v,s}^{(n)} - Y_{v,s}, Y_{w,t}^{(n)})| \leq K_3 |A_n|^{-1} |A_n \setminus A_n^r| \to 0.$$

Thus, by splitting the $v$ sum into $v \in A_n^r$ and $v \in A_n \setminus A_n^r$ we obtain that

$$\limsup_{n \to \infty} |A_n|^{-1} \sum_{v \in A_n} \sum_{w \in A_n} |\mathrm{Cov}(Y_{v,s}^{(n)} - Y_{v,s}, Y_{w,t}^{(n)})| \leq h_2(r)$$

and since $h_2(r) \to 0$ as $r \to \infty$, this shows that

$$\lim_{n \to \infty} |A_n|^{-1} \sum_{v \in A_n} \sum_{w \in A_n} |\mathrm{Cov}(Y_{v,s}^{(n)} - Y_{v,s}, Y_{w,t}^{(n)})| = 0.$$

A similar argument shows that

$$\lim_{n \to \infty} |A_n^{-1}| \sum_{w \in A_n} \sum_{v \in A_n} |\mathrm{Cov}(Y_{v,s}, Y_{w,t}^{(n)} - Y_{w,t})| = 0.$$

Combining the last two limiting expressions, we obtain from (7.9) that

$$|A_n|^{-1} \left( \mathrm{Cov}(S_H^{A_n}(\xi_s^{A_n,\nu}), S_H^{A_n}(\xi_t^{A_n,\nu})) - \mathrm{Cov}(S_H^{A_n}(\xi_s^\nu), S_H^{A_n}(\xi_t^\nu)) \right) \to 0,$$

and combining this with (7.8) we obtain (7.5). ∎

*Proof of Theorem 3.2.* As usual, for any random variable $\xi$ and any $p > 1$ we write $\|\xi\|_p := E[|\xi|^p]^{1/p}$ when the expectation exists. Recall the definition of $Y_z(A, s)$ at (7.1). Define $Y_z'(A, t)$ similarly but with $(\mathcal{P}_\mathbf{0}, \zeta_\mathbf{0})$ resampled (i.e., replaced by an independent copy of $(\mathcal{P}_\mathbf{0}, \zeta_\mathbf{0})$).

For $t > 0$, let $D_t$ denote the set of $v \in \mathbb{Z}^d$ such that there exists $w \in \mathcal{N}_v$ such that $\mathbf{0}$ affects $w$ by time $t$. If we resample $(\mathcal{P}_0, \zeta_\mathbf{0})$ (leaving $((\mathcal{P}_v, \zeta_v), v \in \mathbb{Z}^d \setminus \{0\})$ unchanged), then the set $D_t$ is unchanged; this is because because $(\mathbf{0}, 0)$ is included in $\mathcal{P}$ by (6.6) and any point $(v, T) \in \mathcal{P}$ with $v \neq \mathbf{0}$ that can be reached by a path starting at $(\mathbf{0}, 0)$ can be reached by such a path that does not use any point $(0, T)$ with $T > 0$ (see Section 6). By part (a) of Corollary 6.1, the set $D_t$ is almost surely finite. Also, for any $A \in \mathcal{B}$, if $v \notin D_t$ then we have $Y_v(A, s) = Y'_v(A, s)$ for all $s \leq t$.

We assume there exists $\gamma > 2$ such that (3.3) holds for all $\tau \in I$. Choose $\gamma' \in (2, \gamma)$. For any $A \in \mathcal{B}$ we have by the Minkowski and Hölder inequalities that

$$\left\| \sum_{v \in D_t \cap A} Y_v(A, t) \right\|_{\gamma'} = \left\| \sum_{v \in A} Y_v(A, t) \mathbf{1}_{\{v \in D_t\}} \right\|_{\gamma'} \leq \sum_{v \in \mathbb{Z}^d} \| Y_v(A, t) \mathbf{1}_{\{v \in D_t\}} \|_{\gamma'}$$

$$\leq \sum_{v \in \mathbb{Z}^d} \| Y_v(A, t) \|_\gamma P[v \in D_t]^{(1/\gamma') - (1/\gamma)}$$

and by (3.3) along with Lemma 6.2, this is bounded by a constant that does not depend on $A$.

Let $k \in \mathbb{N}$. Let $(t_1, t_2, \ldots, t_k) \in I^k$ with distinct components, and let $(b_1, \ldots, b_k) \in \mathbb{R}^k$. Consider the linear combination $R_n := \sum_{i=1}^k b_i S_H^{A_n}(\xi_{t_i}^{A_n, \nu})$. Using (3.3) for $\gamma > 2$, we can apply Theorem 3.1 of [47] to deduce that there exists $\sigma^2 \geq 0$, dependent on $(b_i, t_i)_{i=1}^k$, such that as $n \to \infty$,

$$|A_n|^{-1} \text{Var}(R_n) \to \sigma^2, \quad |A_n|^{-1/2}(R_n - ER_n) \xrightarrow{\mathcal{D}} N(0, \sigma^2). \quad (7.11)$$

By the first part of (7.11), and Lemma 7.3,

$$\sigma^2 = \sum_{i=1}^k \sum_{j=1}^k b_i b_j \sigma(t_i, t_j)$$

with $\sigma(s, t)$ given by (3.7). Set

$$\tilde{S}_H^A(\xi_t^{A, \nu}) := S_H^A(\xi_t^{A, \nu}) - E S_H^A(\xi_t^{A, \nu}). \quad (7.12)$$

By the second part of (7.11) and the Cramér-Wold device, the finite-dimensional distributions of the process $(|A_n|^{-1/2} \tilde{S}_H^{A_n}(\xi_t^{A_n, \nu}), t \in I)$ converge to those of a Gaussian process with covariance function $(\sigma(s, t), s, t \in I)$, and this completes the proof. ∎

*Proof of Theorem 3.3.* Again we use the notation $\tilde{S}_H^A(\xi_t^{A, \nu})$ given by (7.12), for the centred version of $S_H^A(\xi_t^{A, \nu})$. We have shown convergence of the finite-dimensional distributions of the sequence of processes $(|A_n|^{-1/2} \tilde{S}_H^{A_n}(\xi_\cdot^{A_n, \nu}), n \geq 1)$, and it remains to prove tightness in $D[0, \tau]$ of this sequence of processes. We follow the standard procedure of demonstrating tightness by estimating moments of increments, as in Billingsley [11], Theorem 15.6.

Let $0 \leq s < t \leq \tau$ with $t - s \leq 1/2$. Changing notation from the previous proof, for $v \in \mathbb{Z}^d$ we set

$$Y_v(s, t) := H(L_v(\xi_t^{A_n, \nu} | \mathcal{N}_v)) - H(L_v(\xi_s^{A_n, \nu} | \mathcal{N}_v)).$$

Then

$$E\left[\left(\tilde{S}_H^{A_n}(\xi_t^{A_n,\nu}) - \tilde{S}_H^{A_n}(\xi_s^{A_n,\nu})\right)^4\right] = E\left[\left(\sum_{v \in A_n}(Y_v(s,t) - EY_v(s,t))\right)^4\right]$$

$$= \sum_{v_1 \in A_n}\sum_{v_2 \in A_n}\sum_{v_3 \in A_n}\sum_{v_4 \in A_n} E\prod_{i=1}^4(Y_{v_i}(s,t) - EY_{v_i}(s,t)). \quad (7.13)$$

For $\mathbf{v} := (v_1, \ldots, v_k) \in (\mathbb{Z}^d)^k$, each $v_i$ ($1 \leq i \leq k$) is an element of $\mathbb{Z}^d$ which we shall refer to as a *d-component* of $\mathbf{v}$. We shall write $v_i \in \mathbf{v}$ as shorthand for '$v_i$ is a $d$-component of $\mathbf{v}$'.

At this point our argument is loosely inspired by the proof of Proposition 8.7 of [18]. Given $\mathbf{v} = (v_1, v_2, \ldots, v_q) \in (\mathbb{Z}^d)^q$ and $r \geq 0$, let us say that $\mathbf{v}$ is *r-connected* if the graph with vertex set $\{1, \ldots, q\}$ and an edge between each pair $(i,j)$ such that $\|v_i - v_j\|_\infty \leq r$ is connected. Define the *gap* of $\mathbf{v}$ to be the smallest $r$ such that $\mathbf{v}$ is $r$-connected. Let $G_r(q,n)$ be the set of $\mathbf{v} \in (A_n)^q$ with gap $r$.

For $\mathbf{v} := (v_1, v_2, v_3, v_4) \in (\mathbb{Z}^d)^4$, and $\mathbf{y} \in (\mathbb{Z}^d)^k$ and $\mathbf{z} \in (\mathbb{Z}^d)^{4-k}$ with $k \in \{1,2,3\}$, we say $\mathbf{y}$ and $\mathbf{z}$ form a *splitting* of $\mathbf{v}$ if there exists nonempty proper subset $J$ of $\{1,2,3,4\}$ such that $\mathbf{y} = (v_i, i \in J)$ and $\mathbf{z} = (v_i, i \in \{1,2,3,4\} \setminus J)$. Define the *distance* between $\mathbf{y}$ and $\mathbf{z}$ to be $\min\{\|y-z\|_\infty\}$, where the minimum is over all $d$-components $y \in \mathbf{y}$ and $z \in \mathbf{z}$.

Suppose $\mathbf{v} = (v_1, \ldots, v_4) \in G_r(4,n)$. Then we assert that $\mathbf{v}$ admits a splitting $(\mathbf{y}, \mathbf{z})$ such that (i) the distance between $\mathbf{y}$ and and $\mathbf{z}$ equals $r$, and (ii) both $\mathbf{y}$ and $\mathbf{z}$ are $r$-connected. To see this, let the edges $\{i,j\}$ of the complete graph on vertex set $\{1,2,3,4\}$ be assigned weights $w_{ij} = \|v_i - v_j\|_\infty$, and let these edges be removed one by one, in order of decreasing weight, until the graph becomes disconnected. Then the last removed edge has weight $r$, and after removing this edge the graph has precisely two components which determine a splitting with properties (i) and (ii) above.

Given $\mathbf{v} = (v_1, \ldots, v_4) \in G_r(4,n)$, let $(\mathbf{y}(\mathbf{v}), \mathbf{z}(\mathbf{v}))$ be a splitting of $\mathbf{v}$ with properties (i) and (ii) of the previous paragraph. Let the choice of $(\mathbf{y}(\mathbf{v}), \mathbf{z}(\mathbf{v}))$, out of all such splittings, be made by an arbitrary deterministic rule. For $\mathbf{y} \in (\mathbb{R}^d)^k$ set

$$\Pi_\mathbf{y} := \prod_{y \in \mathbf{y}}(Y_y(s,t) - EY_y(s,t)). \quad (7.14)$$

Then

$$\sum_{v_1 \in A_n} \sum_{v_2 \in A_n} \sum_{v_3 \in A_n} \sum_{v_4 \in A_n} E \prod_{i=1}^{4}(Y_{v_i}(s,t) - EY_{v_i}(s,t)) = \sum_{r=0}^{\infty} \sum_{\mathbf{v} \in G_r(4,n)} E[\Pi_{\mathbf{v}}]$$

$$= \sum_{r=0}^{\infty} \sum_{\mathbf{v} \in G_r(4,n)} E[\Pi_{\mathbf{y}(\mathbf{v})} \Pi_{\mathbf{z}(\mathbf{v})}]$$

$$= \left( \sum_{r=0}^{\infty} \sum_{\mathbf{v} \in G_r(4,n)} \mathrm{Cov}(\Pi_{\mathbf{y}(\mathbf{v})}, \Pi_{\mathbf{z}(\mathbf{v})}) \right) + \sum_{r=0}^{\infty} \sum_{\mathbf{v} \in G_r(4,n)} E[\Pi_{\mathbf{y}(\mathbf{v})}] E[\Pi_{\mathbf{z}(\mathbf{v})}]. \quad (7.15)$$

If $\mathbf{v} \in G_r(4,n)$, then $\mathbf{y}(\mathbf{v})$ and $\mathbf{z}(\mathbf{v})$ are both $r$-connected. Hence, with $|\cdot|$ denoting cardinality, there exists a constant $K_3$ such that

$$|G_r(k,n)| \leq K_3 |A_n| r^{d-1}(r^d)^{k-2} = K_3 |A_n| r^{d(k-1)-1}. \quad (7.16)$$

Next, we assert that there exists a constant $K_4$ such that for $\mathbf{v} \in G_r(k,n)$,

$$|\mathrm{Cov}(\Pi_{\mathbf{y}(\mathbf{v})}, \Pi_{\mathbf{z}(\mathbf{v})})| \leq K_4 \exp(-r/K_4). \quad (7.17)$$

The proof of this assertion is similar to that of (7.3). Put

$$\hat{\Pi}_{\mathbf{y}} := \prod_{y \in \mathbf{y}} \left( (Y_y(s,t) - EY_y(s,t)) \mathbf{1}_{\{C_{y,t} \subseteq B_{r/3}(y)\}} \right)$$

and choose $\gamma > 6$ such that (3.10) holds. Then by Hölder's inequality and Lemma 6.2,

$$\begin{aligned} E[(\Pi_{\mathbf{y}(\mathbf{v})} - \hat{\Pi}_{\mathbf{y}(\mathbf{v})})^2] &\leq (E[|\Pi_{\mathbf{y}(\mathbf{v})}|^{\gamma}])^{2/\gamma} P[\cup_{y \in \mathbf{y}(\mathbf{v})} \{C_{y,t} \setminus (y + B_{r/3}) \neq \emptyset\}]^{(\gamma-2)/\gamma} \\ &\leq K \exp(-r/K) \end{aligned} \quad (7.18)$$

and likewise for $\mathbf{z}(\mathbf{v})$. Since $\hat{\Pi}_{\mathbf{y}(\mathbf{v})}$ and $\hat{\Pi}_{\mathbf{z}(\mathbf{v})}$ are independent, we obtain

$$|\mathrm{Cov}(\Pi_{\mathbf{y}(\mathbf{v})}, \Pi_{\mathbf{z}(\mathbf{v})})| \leq |\mathrm{Cov}(\Pi_{\mathbf{y}(\mathbf{v})} - \hat{\Pi}_{\mathbf{y}(\mathbf{v})}, \Pi_{\mathbf{z}(\mathbf{v})})| + |\mathrm{Cov}(\hat{\Pi}_{\mathbf{y}(\mathbf{v})}, \Pi_{\mathbf{z}(\mathbf{v})} - \hat{\Pi}_{\mathbf{z}(\mathbf{v})})|$$

and using (7.18), along with the analogous expression for $\mathbf{z}(\mathbf{v})$, and (3.10) and Hölder's inequality gives us (7.17) as asserted.

Next we estimate $\Pi_{\mathbf{y}(\mathbf{v})}$ in terms of $(t-s)$ for small $t-s$. Recall that $\mathbf{v} = (v_1, \ldots, v_4)$ and $\mathbf{y}(\mathbf{v}) = (v_i, i \in J)$ for some nonempty proper subset $J = J(\mathbf{v})$ of $\{1,2,3,4\}$. Expanding the product in (7.14), we obtain

$$\Pi_{\mathbf{y}(\mathbf{v})} = \sum_{J' \subseteq J} \left( \prod_{i \in J \setminus J'} (-E[Y_{v_i}(s,t)]) \right) \times \prod_{j \in J'} Y_{v_j}(s,t), \quad (7.19)$$

where the product over the empty set is 1. Denote by $\tilde{\Pi}_{J'}$ the last factor in (7.19), i.e. $\tilde{\Pi}_{J'} := \prod_{j \in J'} Y_{v_j}(s,t)$. Then $\tilde{\Pi}_{J'} = 0$ unless at least one of the Poisson processes $\mathcal{P}_w, w \in \cup_{i \in J'} \mathcal{N}_{v_i}^+$ has an arrival between time $s$ and time $t$. Hence, if $k$ denotes the number of elements of $J$ then if $J' \neq \emptyset$ there are constants $K_5, K_5'$ such that

$$E[\tilde{\Pi}_{J'}^2] \leq E[|\tilde{\Pi}_{J'}|^{\gamma/k}]^{2k/\gamma} P[\tilde{\Pi}_{J'} \neq 0]^{1-(2k/\gamma)} \leq K_5(t-s)^{1-(2k/\gamma)}$$

so that $E\Pi_{\mathbf{y}(\mathbf{v})}^2 \leq K_5'(t-s)^{1-(2k/\gamma)}$, and similarly $E[\Pi_{\mathbf{z}(\mathbf{v})}^2] \leq K_5'(t-s)^{1-(2(4-k)/\gamma)}$, so that by the Cauchy-Schwarz inequality, there is a constant $K_5''$ such that

$$\mathrm{Cov}(\Pi_{\mathbf{y}(\mathbf{v})}, \Pi_{\mathbf{z}(\mathbf{v})}) \leq K_5''(t-s)^{1-(4/\gamma)}. \tag{7.20}$$

Choose $\alpha > 1/2$ with $\alpha < 1 - (4/\gamma)$ (this is possible because we assume $\gamma > 8$). By (7.16), (7.17) and (7.20), there are constants $K_6$ and $K_7$ such that the first term in the last line of (7.15) satisfies

$$\sum_{r=0}^{\infty} \sum_{\mathbf{v} \in G_r(4,n)} \mathrm{Cov}(\Pi_{\mathbf{y}(\mathbf{v})}, \Pi_{\mathbf{z}(\mathbf{v})})$$

$$\leq K_6 |A_n| \left| \sum_{r \geq 0} r^{3d-1} \min(e^{-r/K_6}, (t-s)^{1-(4/\gamma)}) \right|$$

$$\leq K_6 |A_n| \left( \left( \sum_{0 \leq r \leq (1-(4/\gamma))K_6|\log(t-s)|} r^{3d-1}(t-s)^{1-(4/\gamma)} \right) \right.$$

$$\left. + \sum_{r > (1-(4/\gamma))K_6|\log(t-s)|} r^{3d-1} e^{-r/K_6} \right)$$

$$\leq K_7 |A_n|(t-s)^{\alpha}. \tag{7.21}$$

To estimate the second term in the last line of (7.15), note first that by the same argument as for (7.21), there is a constant $K_8$ such that

$$\sum_{r=0}^{\infty} \sum_{\mathbf{y} \in G_r(2,n)} |E[\Pi_{\mathbf{y}}]| \leq K_8 |A_n|(t-s)^{\alpha}.$$

Then observe that $E[\Pi_{\mathbf{y}(\mathbf{v})}] E[\Pi_{\mathbf{z}(\mathbf{v})}] = 0$ unless both $\mathbf{y}(\mathbf{v})$ and $\mathbf{z}(\mathbf{v})$ have two $d$-components, and hence

$$\sum_{r=0}^{\infty} \sum_{\mathbf{v} \in G_r(4,n)} |E[\Pi_{\mathbf{y}(\mathbf{v})}] E[\Pi_{\mathbf{z}(\mathbf{v})}]| \leq \left( \sum_{r=0}^{\infty} \sum_{\mathbf{y} \in G_r(2,n)} |E[\Pi_{\mathbf{y}}]| \right)^2$$

$$\leq K_8^2 |A_n|^2 (t-s)^{2\alpha}. \tag{7.22}$$

Combining (7.21) and (7.22) in (7.15) we find that for some constant denoted $K_9$, the right hand side of (7.13) is bounded by

$$K_9(|A_n|(t-s)^{\alpha} + |A_n|^2(t-s)^{2\alpha}),$$

and hence by (7.13),

$$E\left[\left(|A_n|^{-1/2}(\tilde{S}_H^{A_n}(\xi_t^{A_n,\nu}) - \tilde{S}_H^{A_n}(\xi_s^{A_n,\nu}))\right)^4\right] \leq K_9\left(\frac{(t-s)^\alpha}{|A_n|} + (t-s)^{2\alpha}\right). \tag{7.23}$$

We now set $U_n(t) := |A_n|^{-1/2}\tilde{S}_H^{A_n}(\xi_t^{A_n,\nu})$. Given $\delta > 0$ with $\delta \leq 1/2$, by (7.23) and the Cauchy-Schwarz inequality we can find $n_0 = n_0(\delta)$ such that for $n \geq n_0$ and $0 \leq t_1 < t \leq t_2 \leq \tau$ with $t_2 - t_1 < \delta$,

$$E[|U_n(t) - U_n(t_1)|^2|U_n(t_2) - U_n(t)|^2] \leq 2K_9\delta^{2\alpha} =: (F(\delta))^{2\alpha}.$$

We slightly modify the argument of proof of Theorem 15.6 of Billingsley [11], using notation $w_F$ from (14.2) of [11] and $w''(U_n,\delta)$ from (14.44) of [11]. Given $\varepsilon > 0$ and $\delta > 0$, choose $n_0(\delta)$ as above. By following pages 129–130 of [11] we obtain for $n \geq n_0$ that (cf. (15.30) of [11])

$$P[w''(U_n,\delta) \geq \varepsilon] \leq \frac{2K'}{\varepsilon^4}(F(\tau) - F(0))(w_F(2\delta))^{2\alpha-1} \tag{7.24}$$

where the constant $K'$ does not depend on $\delta$. Given $\varepsilon > 0$ and $\eta > 0$, we can choose $\delta$ so that

$$(2K'/\varepsilon^4)(F(\tau) - F(0))(w_F(2\delta))^{2\alpha-1} < \eta,$$

so by (7.24), we have (15.22) of [11] for $n \geq n_0(\delta)$, and we can then follow [11] to get the desired tightness and convergence in $D[0,\tau]$.

With the limiting Gaussian process denoted $(Z(t), 0 \leq t \leq \tau)$, we have from (7.23) and Theorem 3.2, along with the Skorohod coupling [38] and Fatou's lemma, that $E[(Z(t) - Z(s))^4]$ is bounded by $K_9(t-s)^{2\alpha}$, so that by the Kolmogorov-Chentsov criterion (see e.g. [38]), the process $G(\cdot)$ admits a version with continuous sample paths. ∎

**Acknowledgements.** Some of this work arose following discussions with Joseph Yukich. Travel for these discussions was partly funded by the London Mathematical Society.

We thank the referee for pointing out a flaw in the presentation of Theorem 2.1 in an earlier version of this paper. We thank Sana Louhichi and Xin Qi for providing preprint versions of their papers [18] and [54] repsectively.